\documentclass[11pt,reqno]{amsart}
\usepackage {amsmath,amsfonts,enumerate, epsfig, epsf,  indentfirst,  graphicx,euscript, amssymb,amsthm,amscd}

\usepackage[latin1]{inputenc}

\makeindex

\theoremstyle{plain}
\newtheorem{proposition}{Proposition}

\newtheorem{remark}{Remark}
\theoremstyle{definition}
\newtheorem{definition}{Definition}
\newtheorem{theorem}{Theorem}

\newtheorem{problem}{Problem}

 \def \ue{\"u\^e}
 \def\b{\mathbf}

\def\uu{(u_1,u_2,u_3)}

\def\oa#1#2{\omega_{#1#2}}

 \def\ue#1#2{u^{#1}_{#2}}
\def\ui#1{u_{#1}}

\def\we{\wedge}

\def\lp#1{${\mathcal L}_{#1}(\alpha)$}
\def\fp#1{${\mathcal F}_{#1}(\alpha)$}

\def\prf#1#2{\dfrac{\partial #1}{\partial #2}}

\def\re{\mathbb R}

\def\lise{partially umbilic  curve{}}
\vspace{-.9cm}
\title[Recent Developments on Lines of  Curvature]{Lines of
 Curvature  on Surfaces, Historical Comments and Recent Developments}

\author[J. Sotomayor and R. Garcia ]{Jorge Sotomayor and Ronaldo  Garcia}


 \thanks{The  authors are fellows of CNPq and
done this work under the project CNPq 473747/2006-5}

\begin{document}

  \begin{abstract}
This survey  starts with
the historical landmarks leading to the study
of  principal configurations on surfaces, their structural stability
and further  generalizations. Here it is  pointed out that in the
work of Monge, 1796,  are found elements of the qualitative theory
of differential equations ({\it QTDE}), founded by Poincar\'e in 1881.
Here are also outlined a number of recent results developed after the assimilation into the subject  of
 concepts and problems from the {\it QTDE} and Dynamical Systems, such as  Structural Stability,
  Bifurcations and Genericity, among others,  as well as extensions to higher dimensions.
  References to original works are given and
 open problems are proposed at the end of some sections.
\end{abstract}

\vspace{-1.66cm}
 \maketitle
 \vspace{-1.2cm}
 {\Small\tableofcontents}

\section{ Introduction }
\label{sec:1}

The book on differential geometry of  D. Struik   \cite {St}, remarkable for
its historical notes,
contains key  references to the classical works on {\em principal
curvature lines} and their {\em umbilic } singularities due  to L.
Euler \cite{le}, G. Monge \cite{mo}, C. Dupin \cite{du}, G. Darboux \cite{da} and A. Gullstrand \cite{gul}, among
others (see  \cite{Sp} and, for additional references, also
\cite{r}). These papers
---notably that of Monge, complemented with Dupin's--- can be
connected with aspects of the {\em  qualitative theory of
differential equations} ({\em QTDE} for short) initiated by H.
Poincar\'e \cite{pome} and consolidated with the study of the {\it
structural stability and genericity} of differential equations in
the plane and on surfaces, which was made systematic from $1937$
to $1962$ due to
 the seminal works of Andronov  Pontrjagin and Peixoto (see \cite{ap} and \cite{mp}).

This survey  paper
 discusses the  historical sources for the  work on
 the  structural stability of  principal curvature
 lines and umbilic points,  developed by
  C. Gutierrez and J. Sotomayor  \cite{gs1, gsln, gs2}.
  Also it addresses other kinds of geometric foliations studied
      by  R. Garcia and
    J. Sotomayor \cite{gsalg, ax, me}. See   also  the papers devoted to
    other  differential
   equations
   of classical geometry: the  asymptotic lines \cite {gaguso, gas}, and  the  arithmetic, geometric
    and harmonic  mean curvature lines \cite{m, g, h, me}.

In the historical comments posted in  \cite{arXho} it is  pointed out that in the work of Monge,
\cite{mo}, are found elements of the {\em QTDE}, founded by Poincar\'e in \cite{pome}.

The present  paper contains a reformulation of
 the essential historical aspects  of  \cite {moe, moe2}. Following the thread of Structural Stability and the {\em QTDE}
 it discusses pertinent extensions and    updates  references.  At the end of  some sections   related open problems are
 proposed and commented.

 Extensions of the  results outlined in section  \ref{sec:gs} to
surfaces with generic critical points, algebraic surfaces, surfaces and  $3$  dimensional manifolds
in ${\mathbb R} ^4$  have been
achieved recently (see, for example, \cite{ga2, gaguso, gms, gsalg, gas, ax,  gsjap, gsc2,
 mello}).  An  account of these recent developments  will be given here.

\section{Historical Landmarks}\label{sec:hl}
\subsection{ The Landmarks before Poincar\'e: Euler, Monge and  Dupin}
\label{sec:2}\hskip 1cm

\begin{figure}[htb]
\begin{center}
\includegraphics[height=5cm, height=3.5cm]{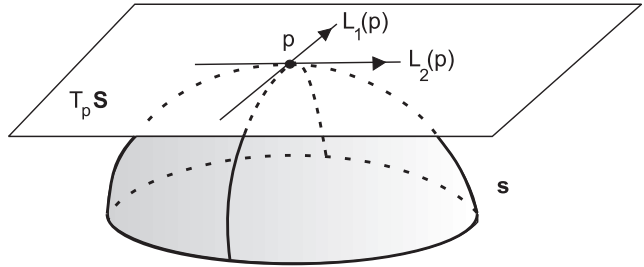}
\end{center}
\caption{Principal Directions. \label{fig:1}}
\end{figure}

\medskip
\noindent {\bf Leonhard Euler} (1707-1783) \cite{le}, founded
the curvature  theory of surfaces. He defined the {\em normal
curvature} $k_n (p,L)$ on an oriented surface {\bf S} in a tangent
direction $L$ at a point $p$ as the
 curvature, at $p$,   of the planar curve of intersection of the surface with the
 plane generated by the line $L$ and the positive unit normal $N$ to the surface at $p$.
 The  {\em   principal curvatures} at $p$ are the extremal values of  $k_n (p,L)$ when
 $L$ ranges over  the tangent directions through $p$. Thus,     $k_1(p)  = k_n (p, L_1) $
 is the {\em  minimal}   and $k_2 (p)= k_n (p, L_2) $ is the {\em maximal normal curvatures},  attained
 along the  {\em principal directions}: $L_1 (p)$,  {\em the minimal}, and  $L_2 (p)$,
 {\em  the maximal} (see Fig. \ref{fig:1}).

Euler's formula expresses the normal curvature $k_n (\theta)$
along a direction making angle $\theta$ with the minimal principal
direction  $L_1$ as  $k_n (\theta)=k_1 \cos ^2 \theta  + k_2
\sin^2 \theta$.

Euler, however,  seems to have not  considered the integral curves
of the principal line fields $L_i : p\to L_i (p), \,  i\, =1,\, 2,
\, $ and overlooked the role  of the {\em umbilic points}  at
which the principal curvatures coincide and the line fields are
undefined.

\vskip 0.2cm

\noindent {\bf Gaspard Monge} (1746-1818)
 found the family of integral curves of the {\em principal line
fields}  $L_i ,    i\, =1,\, 2, \,  $ for the case of the triaxial
ellipsoid
$$\frac{x^2}{a^2} + \frac{y^2}{b^2} + \frac{z^2}{c^2}  -
1 = 0, \,\,\, a>b>c>0.$$
In doing this, by direct integration of
the differential equations of the principal curvature lines, circa
1779, Monge was led to the first example of a  {\em foliation with
singularities} on a surface  which (from now on) will be called
the {\em principal configuration} of an  oriented surface.
The singularities consist  on the {\em umbilic points},
mathematical term he coined to designate those at which
the principal curvatures coincide and the line fields are undefined.

\begin{figure}[htb]
\begin{center}
\includegraphics[height=4cm]{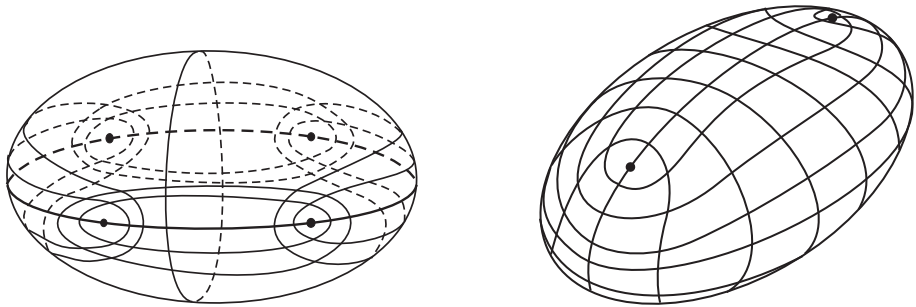}
\end{center}
\caption{Monge's Ellipsoid, two perspectives. \label{fig:2}}
\end{figure}

The
Ellipsoid, endowed with its principal configuration, will be
called Monge's Ellipsoid (see Fig. \ref{fig:2}).

\vskip 0.2cm

The motivation found in Monge's paper \cite {mo} is a complex
interaction of esthetic and  practical considerations and of the
explicit desire  to apply the results of his mathematical research
to real world problems. The principal configuration on the
triaxial ellipsoid appears in Monge's proposal  for the dome of
the Legislative Palace for the government of the French
Revolution, to be built over an elliptical terrain. The lines of
curvature are the guiding  curves for the workers to put the
stones. The umbilic points, from which were to hang the candle
lights, would also be the reference points below which to put the
podiums for  the speakers.

\vskip 0.2cm

\noindent {\bf Commenting  Monge's work under the perspective of the {\em QTDE}}
 The ellipsoid depicted in Fig. \ref{fig:2} contains some of the  typical
features of the qualitative theory of differential equations
discussed briefly in {\bf a)} to {\bf d)} below:

{\bf a)  Singular Points and Separatrices.}   The umbilic points
play the role of singular points for the principal foliations,
each of them has one separatrix for each principal foliation. This
separatrix produces a connection with another umbilic point of the
same type, for which it is also a separatrix, in fact  an umbilic
separatrix connection.

\vskip 0.2cm

{\bf b)  Cycles.} The configuration  has principal cycles. In
fact, all the principal lines, with the exception of the four umbilic
connections, are periodic. The cycles fill a cylinder or annulus,
for each foliation. This pattern is common to all classical
examples, where no surface exhibiting an isolated cycle was known.
This fact seems to be derived from the symmetry of the surfaces
considered, or from the integrability that is present in the
application of Dupin's Theorem for triply orthogonal families of
surfaces.

 As was  shown
in \cite{gs1}, these configurations are
exceptional; the generic
principal cycle  for a smooth surface is  a  hyperbolic limit
cycle (see below).

\vskip 0.2cm
 {\bf c) Structural Stability (relative to quadric
surfaces).} The principal configuration remains qualitatively
unchanged under small perturbations on the coefficients of the
quadratic polynomial that defines the surface.
 \vskip 0.2cm

 {\bf d) Bifurcations.}
The drastic changes in the principal configuration exhibited by
the transitions from a sphere, to an ellipsoid of revolution and
to a triaxial ellipsoid (as in  Fig. \ref{fig:2}), which after  a
very small perturbation, is a simple form of a bifurcation
phenomenon.

\vskip 0.2cm

\begin{figure}[htb]
\begin{center}
\includegraphics[height=5cm ]{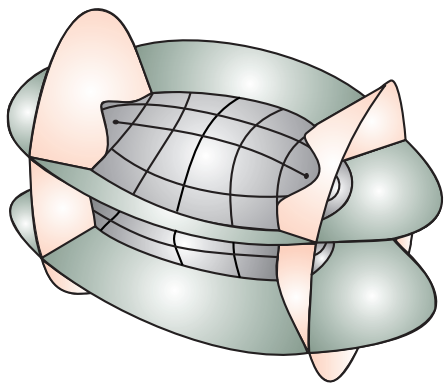}
\end{center}
\caption{Dupin's Theorem. \label{fig:3}}
\end{figure}

\noindent {\bf Charles Dupin } (1784-1873) considered the surfaces
that  belong to {\em triply orthogonal surfaces,}  thus extending
considerably  those whose  principal configurations can be found
by integration. Monge's Ellipsoid belongs to the family of {\em
homofocal quadrics} (see \cite{St} and  Fig.  \ref{fig:3}).

The conjunction of  Monge's analysis and
Dupin extension  provides the first global theory of integrable
principal configurations, which for quadric surfaces gives those
which are also {\em principally structurally stable } under small
perturbations of the coefficients of their quadratic  defining
equations.

\begin{theorem}\label{th:q2} \cite {moe}
In the space of oriented quadrics, identified with the
nine-dimensional sphere,  those having principal structurally
stable configurations are open and dense.
\end{theorem}

\vskip 0.2cm

\noindent{\bf Historical Thesis in \cite{arXho}}.  {\em The global study of  lines
of principal curvature leading to Monge's Ellipsoid, which is
analogous of the phase portrait of a differential equation,
contains elements of Poincar\'e's  QTDE, 85 years earlier.}

\vskip 0.2cm
 This connection  seems to have been  overlooked by Monge's
scientific
 historian Ren\' e  Taton (1915-2004) in his remarkable book
\cite{ta}.

\subsection { Poincar\'e  and  Darboux }\label{sc:3}

The exponential role played
by {\bf Henri Poin\-ca\-r\'e} (1854-1912) for the {\em QTDE} as well
as for other branches of mathematics is well known and has been
discussed and analyzed in several places (see for instance
\cite{fb} and \cite{ppr}).

\begin{figure}[htb]
 \begin{center}
 \includegraphics[height=2.5cm ]{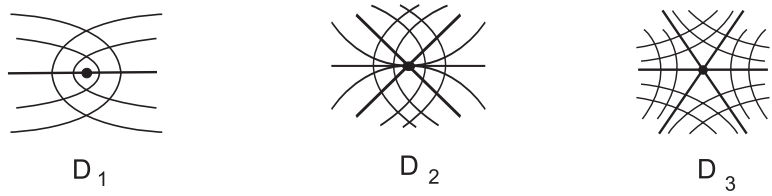}
 \caption{ Darbouxian Umbilics. \label{fig:4}}
   \end{center}
 \end{figure}

\begin{figure}[htbp]
 \begin{center}
 \includegraphics[height=4cm ]{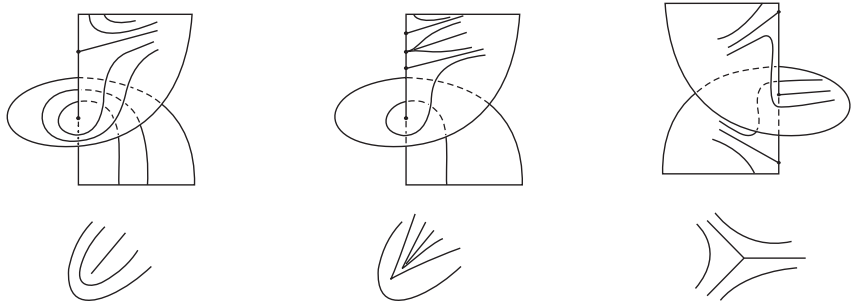}
 \caption{Lie-Cartan Resolution of Darbouxian Umbilics. \label{fig:5}}
   \end{center}
 \end{figure}

Here we are concerned with his M\'emoires  \cite{pome}, where he
laid the foundations of the {\em QTDE}.  In this work Poincar\'e
determined the form of the solutions of planar analytic
differential equations near their  {\em foci},  {\em nodes} and
{\em saddles}. He also  studied   properties of the solutions
around cycles and, in the case of polynomial differential
equations, also the behavior at infinity.
 \vskip 0.2cm \noindent
{\bf Gaston Darboux} (1842-1917) determined the  structure of the
lines of principal curvature near a {\em generic} umbilic point.
In his note \cite{da}, Darboux uses the theory of singularities of
Poincar\'e.  In fact, the
Darbouxian umbilics are those whose resolution
by {\em blowing up} reduce only to
 saddles and nodes (see Figs. \ref{fig:4} and
\ref{fig:5}).

\vskip 0.2cm

Let $\b  p_0\in \b S$ be an umbilic point. Consider a chart
$(u,v):(\b S,\b  p_0)\to ({\mathbb R}^2,\b  0)$ around it, on which the
surface has the form of the graph of a function such as
$$
\frac{k}{2}(u^2+v^2)+ \frac{a}{6}u^3+ \frac{b}{2}uv^2+
\frac{c}{6}v^3+ O[(u^2+v^2)^2].
$$

\noindent This is achieved by projecting $\b S$ onto $T_{\b p_0}\b S $ along $N(\b  p_0)$ and choosing there an orthonormal chart
$(u,v)$ on which the coefficient of the cubic term $u^2v$
vanishes.

An umbilic point is called {\it Darbouxian\/}  if, in the above
expression, the following 2 conditions {\bf T)} and {\bf D)} hold:

{\bf T)} $ b(b-a)\neq0 $,

{\bf D)} either

 {\rm D}$_1$: $a/b > (c/2b)^2+2$,

 {\rm D}$_2$: $(c/2b)^2+2 >a/b >1$,  $a\neq2b$,

 {\rm D}$_3$: $a/b <1. $

The corroboration  of the pictures in Fig.  \ref{fig:4}, which
illustrate the principal configurations near Darbouxian umbilics,
has been given in \cite{gs1, gs2}; see also \cite{bf} and  Fig.
\ref{fig:5} for the Lie-Cartan resolution of a Darbouxian umbilic.

The subscripts  refer to the number of {\em umbilic separatrices,}
 which are the curves, drawn with heavy lines, tending to the
umbilic point and separating regions whose principal lines have
different patterns of approach.

\vskip 0.2cm

\subsection{ Principal Configurations on Smooth Surfaces in  $\mathbb R^3$}\label{sec:gs}

\hskip 2cm
\medskip

After the seminal work of Andronov-Pontrjagin \cite{ap} on
structural stability of differential equations in the plane and
its extension to  surfaces by Peixoto \cite{mp} and in view of
the discussion on Monge's Ellipsoid formulated above, an  inquiry
into  the characterization of the oriented surfaces  $\mathbf S$
whose principal configuration are structurally stable under small
$C^r$ perturbations, for $r\geq3$, seems unavoidable.

\vskip 0.2cm

 Call $\Sigma(a,b,c,d)$ the set of smooth compact oriented
surfaces $\b S$ which verify the following  conditions.

{\bf a)} All umbilic points are Darbouxian.
\vskip 0.1cm
 {\bf b)} All principal cycles are hyperbolic.
 This means that the corresponding   return map is
{\it hyperbolic\/}; that is,  its derivative is different from 1.
It has been shown in \cite{gs1} that hyperbolicity of a principal
cycle $\gamma$ is equivalent to the requirement  that
\vspace{-0.3cm}
$$
\int_\gamma \frac{d\mathcal H} {\sqrt{\mathcal H^2-\mathcal K}}
\neq 0,
$$
\noindent where $\mathcal H = (k_1 +k_2)/2$ is the mean curvature
and $\mathcal K = k_1 k_2$ is the gaussian curvature.

\vskip 0.1cm
 {\bf c)} The limit set of every principal line is
contained in the set of umbilic points and principal cycles of $\b
S$.

The $\alpha$-(resp. $\omega$) {\it limit set\/}  of an oriented
principal line $\gamma$, defined on its maximal interval $\mathcal
I=(w_-,w_+)$ where it is parametrized by arc length $s$, is the
collection $\alpha(\gamma)$-(resp. $\omega(\gamma)$) of limit
point sequences of the form $\gamma(s_n)$, convergent in $\b S$,
with $s_n$ tending to the left (resp. right) extreme of $\mathcal
I$. The {\it limit set\/}  of $\gamma$ is the set
$\Omega=\alpha(\gamma)\cup\omega(\gamma)$.

 Examples of surfaces with {\it non
trivial recurrent\/}  principal lines, which violate condition
{\bf c} are  given in \cite{gsln, gs2} for ellipsoidal and toroidal surfaces.

 There are no examples of
these situations in the classical geometry litera\-tu\-re.

\vskip 0.1cm
 {\bf d)} All umbilic separatrices are separatrices of
a single umbilic point.

 Separatrices which violate {\bf d} are called {\it umbilic
connections\/}; an example can be seen in the ellipsoid of Fig.
\ref{fig:2}.

\vskip 0.2cm To make precise the  formulation  of the next
theorems, some topological notions must be defined.

 A sequence $\b S_n$ of
surfaces {\it converges in the $C^r$ sense} to a surface $\b S$
provided there is a sequence of real functions $  f_n$ on $\b S$,
such that $\b S_n= (I+ f_n N)(\b S)$ and $f_n$ tends to $0$ in the
$C^r$ sense; that is, for every chart $(u,v)$ with inverse
parametrization $X$, $f_n\circ X$ converges to $0$, together with
the partial derivatives of order $r$, uniformly on compact parts
of the domain of  $X$.

A set $\Sigma$ of surfaces is said to be {\it open\/}  in the
$C^r$ sense if every sequence $\b S_n$ converging to $\b S$ in
$\Sigma$ in the $C^r$ sense is, for $n$ large enough, contained in
$ \Sigma$.

A set $\Sigma$ of surfaces is said to be {\it dense}  in the $C^r$
sense if, for every surface $\b S$, there is a sequence $\b S_n$
in $\Sigma$ converging to $\b S$ the $C^r$ sense.

A surface $\b S$ is said to be $C^r$-{\it principal structurally
stable\/}  if for every sequence $\b S_n$ converging to $\b S$ in
the $C^r$ sense, there is a sequence of homeomorphisms $H_n$ from
$\b S_n$ onto $\b S$, which converges to the identity of $\b S$,
such that, for $n$ big enough, $ H_n$ is  a principal equivalence
from $\b S_n$ onto $\b S$. That is, it maps $\b U_n$, the umbilic
set of $\b S_n$, onto $\b U$, the umbilic set of $\b S$, and maps
the lines of the principal foliations $\b F_{i,n}$, of $\b S_n$,
onto those of $\b F_i$, $i=1,2$, principal foliations for $\b S$.

\begin {theorem} \label{th:2} (Structural Stability of Principal
Configurations \cite{gs1, gs2})
 The set of surfaces
$\Sigma(a,b,c,d)$ is open in the $C^3$ sense and each of its
elements is $C^3$-principal structurally stable.
\end{theorem}

\begin {theorem} \label{th:3} (Density of Principal Structurally Stable
Surfaces, \cite{gsln, gs2})
 The set $\Sigma(a,b,c,d)$ is dense in the $C^2$ sense.
\end{theorem}

\vskip 0.2cm

To conclude this section two open problems are proposed.

\begin{problem} \label{pro:1}
Raise from 2 to 3 the differentiability  class in the density
Theorem \ref{th:3}.
\end{problem}

This remains among the most intractable questions in this subject,
involving  difficulties of Closing Lemma type, \cite{pu}, which
also permeate other differential equations of classical geometry,
\cite{me}.

\begin{problem} \label{pro:3toro}
Is it possible to have a smooth  embedding  of the Torus $\mathbb T^2$ into
$\mathbb R^3$ with  both  maximal and minimal non-trivial
recurrent principal curvature  lines?

\end{problem}

Examples with either maximal or minimal principal recurrences can
be found in \cite{gsln, gs2}.

\section{ Curvature Lines near Umbilic Points}\label{sec:lcumb}

  The purpose  of this section is to present the simplest qualitative
changes --bifurcations-- exhibited by  the principal
configurations under small perturbations of an immersion which
violates   the  Darbouxian structural stability
condition on umbilic points.

It will be presented the two codimension one umbilic points {\rm D}$^{1}_{2}$  and {\rm D}$^{1}_{2,3}$, illustrated in
Fig.  \ref{fig:semidarboux} and  the four codimension two umbilic points {\rm D}$^{2}_{1}\,$, {\rm D}$^{2}_{2p} \,$, {\rm D}$^{2}_{3}\,$ and {\rm D}$^{2}_{2h}$, illustrated in Figs. \ref{fig:4} and \ref{fig:d13}.

The superscript stands for the codimension which is the
minimal  number of parameters on which depend the families of
immersions exhibiting
 persistently the pattern.    The subscripts stand for the number of
  separatrices approaching the umbilic. In the first case, this number
   is the same for both the minimal and maximal principal curvature
   foliations. In the second case, they are not equal and, in our notation,
   appear  separated by a comma.
     The symbols $p$,
for {\it parabolic},  and $h$, for {\it hyperbolic}, have been
added to the subscripts above in order to distinguish types that
are not discriminated only by the number of separatrices.

\subsection{ Preliminaries on   Umbilic Points} \label{sec:2u}\hskip 2cm
\medskip

The following assumptions will hold from now on.

 Let $p_0$ be an
umbilic point of an immersion $\alpha$ of an oriented surface
$\mathbb M $ into $\mathbb R ^3$, with a once for all fixed
orientation. It will be assumed that $\alpha$ is of class $C^k ,\;
k \geq \, 6$.
 In a local Monge  chart
 near $p_0$  and a positive  adapted $3-$frame,
  $\alpha$ is given by $\alpha (u,v)=(u,v,h(u,v))$,
where

 \begin{equation}\label{eq:m1}
  \aligned h (u,v) &= \frac k2 (u^2+v^2) + \frac a6 u^3 +\frac
b2 u v^2 +\frac c6 v^3
+\frac A{24} u^4 + \frac B6 u^3 v\\
 &+\frac C4 u^2v^2 + \frac D6 u v^3 + \frac E{24} v^4
+\frac{a_{50}}{120}u^5+ \frac{a_{41}}{24}u^4v\\
&+ \frac{a_{32}}{12}u^3v^2 +\frac{a_{23}}{12}u^2v^3+
\frac{a_{14}}{24}uv^4+\frac{a_{05}}{120}v^5+ h.o.t
\endaligned
\end{equation}

Notice that, without loss of generality,  the term $u^2 v$ has
been eliminated from this expression by means of a rotation in the
$(u,v)-$ frame.

 According to \cite{Sp} and \cite{St}, the differential equation of lines of
curvature in terms of  $I=Edu^2 +2Fdudv+Gdv^2$ and  $II=edu^2
+2fdudv+gdv^2$ around $p_0$ is:

\begin{equation}\label{eq:lc}
(Fg-Gf)dv^2+(Eg-Ge)dudv+(Ef-Fe)du^2=0.
\end{equation}

\noindent  Therefore the  functions $L=Fg-Gf$, $M=Eg-Ge$ and
$N=Ef-Fe$ are:
$$\aligned L &= h_uh_v h_{vv}- (1+h_v^2)h_{uv}\\
M &= (1+h_u^2)h_{vv}- (1+h_v^2)h_{uu}\\
N &= (1+h_u^2)h_{uv} - h_u h_v h_{uu}.\endaligned $$

  Calculation gives
{\small
\begin{equation}\label{eq:lmn}
\aligned L &= -b v -\frac 12 B
u^2 - (C-k^3)u v -\frac 12 D v^2
-\frac{a_{41}}6u^3\\
&+\frac 12(4bk^2+k^2a-a_{32})u^2v
+\frac 12(3k^2c-a_{23})uv^2 -\frac 16(a_{14}+3bk^2)v^3
 + h.o.t\\
 M &= (b-a)u + c v
+ \frac 12(C-A+2k^3)u^2 + (D-B) u v+ \frac 12 ({E-C} -2k^3)v^2
\\ & +\frac 16[6bk^2(a+b) +a_{32}-a_{50}]u^3+
\frac 12(a_{23}+2ck^2-a_{41})u^2v\\
&+\frac 12[a_{14}-a_{32}-2 k^2(a+b)]uv^2+
\frac 16(a_{05}-a_{23}-6ck^2)v^3+
 h.o.t\\
N&= b v + \frac 12 B  u^2 + (C-k^3)u v + \frac 12 D  v^2
+ \frac 16  a_{41} u^3\\
&+\frac 12(a_{32}-3ak^2)u^2v
+\frac 12(a_{23}-k^2c)uv^2
+\frac 16(a_{14}-3bk^2)v^3
+ h.o.t \endaligned \end{equation}
}

\subsection  {Umbilic Points of Codimension One} \label{ssec:cod1}

\hskip 2cm
\medskip

A characterization of umbilic points of codimension one is the following
theorem  announced in \cite{gs-5} and proved in \cite{sell}.

\begin{figure}[htbp]
 \begin{center}
  \includegraphics[angle=0, height=4cm]{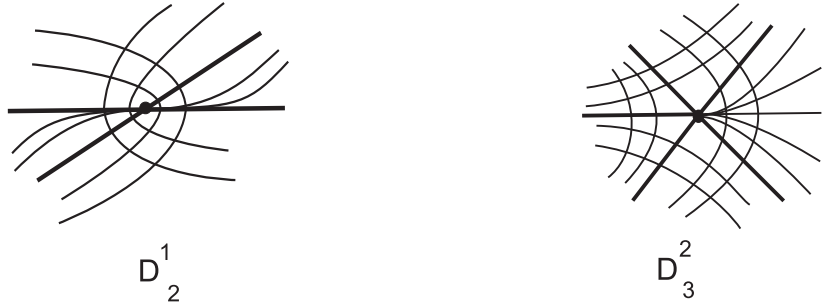}\hskip 1cm
  \caption{\label{fig:semidarboux} Principal curvature lines near the umbilic points
   {\rm D}$^{1}_{2}$, left, and
  {\rm D}$^1_{2,3}$, right,
 and their separatrices. }
 \end{center}
  \end{figure}

\begin{theorem}  {\cite{sell, gs-5}}\label{th:cod1}
Let $p_0$ be an umbilic point and consider
$\alpha(u,v)=(u,v,h(u,v))$
 as in equation  \eqref{eq:m1}.
 Suppose the following conditions hold:

\begin{itemize}
   \item[{{\rm D}$^{1}_{2}$})]\quad $ cb(b-a)\ne 0\,$  and  either
\quad $ (\frac{c}{2b})^2 - \frac ab +1 = 0\quad$   or \quad $
a=2b.$
 \item[{{\rm D}$^{1}_{2,3}$})] \quad $b=a\ne 0\quad \text{and} \quad
   \chi={cB}  -({C-A}  + 2k^3)b\ne 0. $
\end{itemize}

 Then the behavior of lines of
curvature near the  umbilic point $p_0$ in  cases {\rm D}$^{1}_{2}$ and
{\rm D}$^{1}_{2,3}$,
 is as illustrated  in Fig.  \ref{fig:semidarboux}.

\end{theorem}

\vspace{-.3cm}

\begin{figure}[htbp]
  \begin{center}
  \hskip 1cm
  \includegraphics[angle=0, height=5.0cm, width=10cm]{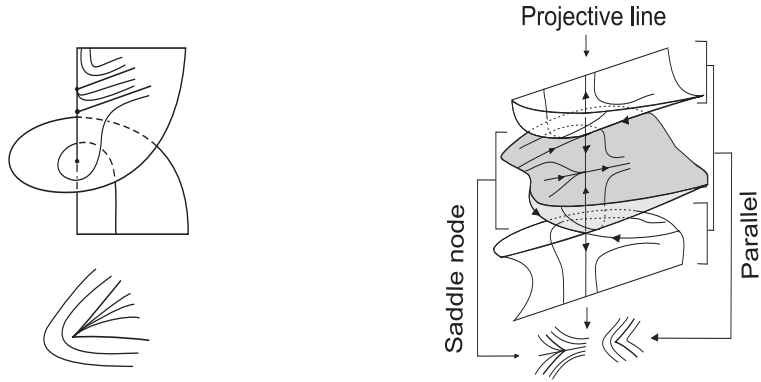}
  \caption{Lie-Cartan suspension   {\rm D}$^{1}_{2}$, left, and
   {\rm D}$^{1}_{2,3}$, right.
  \label{fig:lcd}}
    \end{center}
  \end{figure}

In  Fig.  \ref{fig:lcd} is illustrated the behavior of the Lie-Cartan resolution of
  the semi Darbouxian umbilic points.

     Global effects, due to
umbilic bifurcations, on these configurations such as the
appearance and annihilation of
  principal cycles  were studied in \cite{sell} and
  are outlined in section \ref{sec:ucloop}.
\medskip

\subsection  { Umbilic Points of Codimension Two} \label{ssec:cod2}
\hskip 2cm
\medskip

The characterization of umbilic points of codimension two, generic in biparametric families of immersions, were established in \cite{gsc2} and will be reviewed in  Theorem \ref{th:ud2} below.

\begin{figure}[htbp]
       \begin{center}
 \includegraphics[angle=0, height= 4.0cm ]{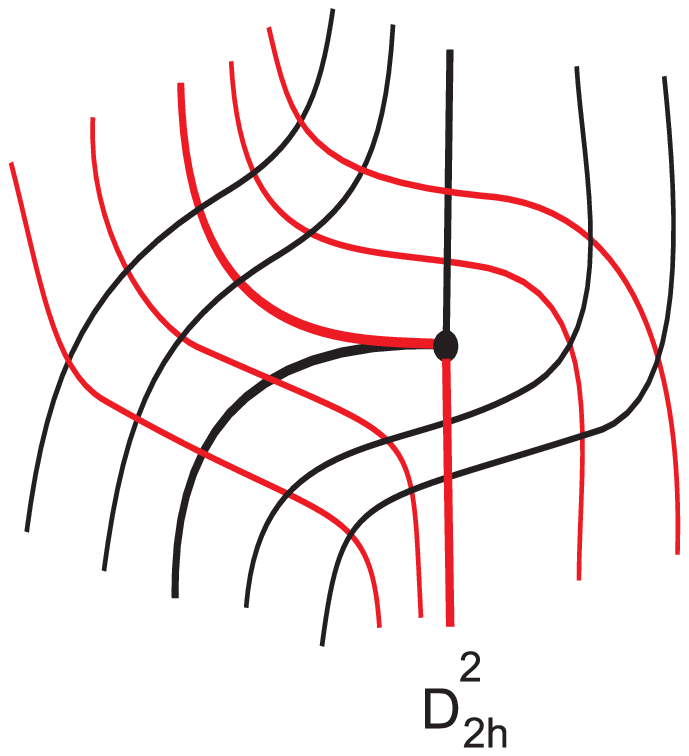}\hskip 2cm
 \includegraphics[height=5.0cm ]{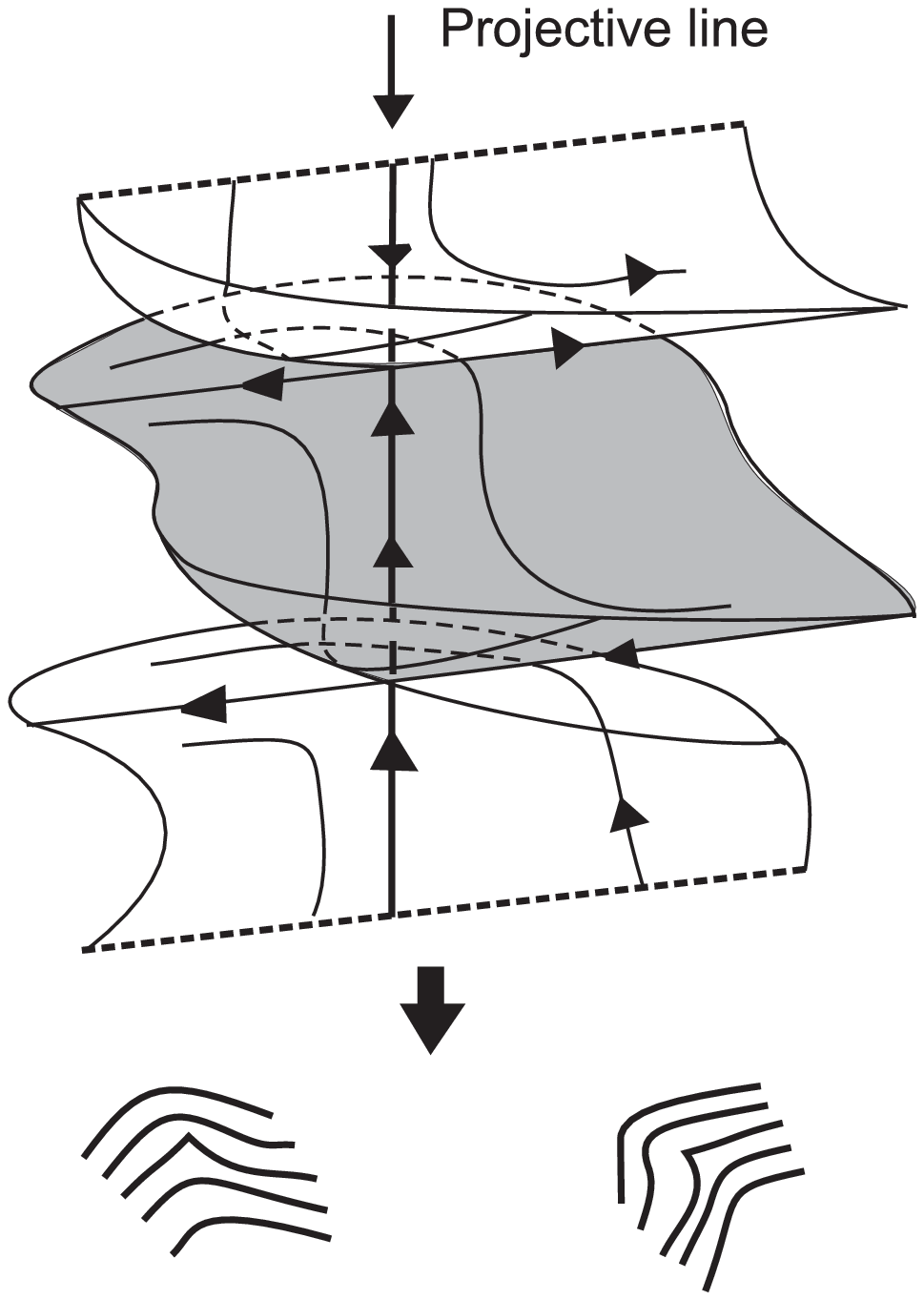}
  \caption{Lines of curvature near the umbilic
   point  {\rm D}$^{2}_{2h}$, left,  and associated Lie-Cartan suspension, right. \label{fig:d13}}
    \end{center}
  \end{figure}

\begin{theorem} \cite{gsc2} \label{th:ud2} Let $p_0$ be an umbilic point
 and $\alpha(u,v)=(u,v,h(u,v))$
 as in equation \eqref{eq:m1}.
\begin{itemize}
\item[a)]
Case  {\rm D}$^{2}_{1}$):  If  $\;\;c=0$\;\; and\;\; $a=2b\ne 0$,

\noindent then  the configuration  of  principal lines  near $p_0$
is topologically equivalent to    that of  a  Darbouxian  {\rm D}$_1$
 umbilic point. See  Fig.  \ref{fig:4}, left.

\item[b)] If $a=b\ne 0$,
$ \chi = cB- b(C-A +2 k^3)=0 $
 and

\noindent $\xi=  12k^2 b^3+ (a_{32}-a_{50})b^2+(3B^2-3BD
-ca_{41})b+ 3 cB(C-k^3)  \ne 0, \, $

\noindent then the principal configurations of lines of curvature
fall into one of the  two cases:
\begin{itemize}
\item[i)]
 Case {\rm D}$^{2}_{2p}$: $\xi b<0 $, which is topologically a  {\rm D}$_2$  umbilic and
\item[ii)] Case  {\rm D}$^{2}_{3}$:\; $\xi b>0$, which is topologically a
  {\rm D}$_3$  umbilic.
\end{itemize}
See Fig.  \ref{fig:4},  center and right, respectively.

\item[c)]
Case {\rm D}$^{2}_{2h}$: If $a=b=0$ and $cB\ne 0$ or if  $b=c=0$ and $aD\ne 0$, then the
principal configuration near $p_0$ is as in  Fig.  \ref{fig:d13}.

\end{itemize}
\end{theorem}

\vspace{-.5cm}

\subsection{ Umbilic Points of Immersions with Constant Mean Curvature } \label{sc:ucmc}  \hskip 2cm

\begin{figure}[ht]
\begin{center}
\includegraphics[angle=0, height=2.7cm]{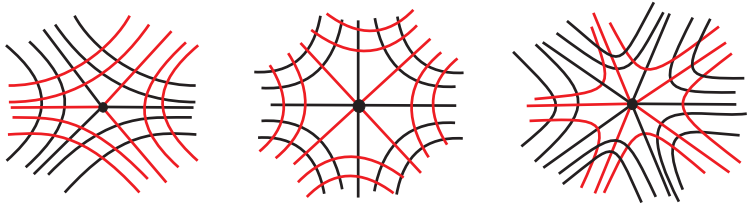}
\end{center}
\caption{ \label{fig:umbh} Curvature lines near   isolated  umbilic points
   of immersions with constant mean curvature.}
\end{figure}

\medskip
\vspace{-.3cm}
Let  $(u,v)$  be isothermic coordinates in a neighborhood of an isolated umbilic point $p=0$ of an immersion $\alpha:\mathbb M^2\to \mathbb R^3$
with constant mean curvature ${\mathcal H}_\alpha$.  In terms of $\phi= (e-g)/2-if $, and $w = u + iv$, the equation  of principal lines is written
as  $\;Im(\phi(w)dw^2)=0.$ See also \cite{gstr} and \cite{hopf}.
The index of an isolated umbilical point with complex coordinate
$w = 0 $ is equal to $-n/2$, where $n$ is
the order of the zero of $ \phi $ at $w = 0$.
 There are $ n + 2 $ rays $ L_0,\; L_1,\ldots,L_{n+1}$
through $0 \in T_p\mathbb M$, of which two consecutive rays make an angle of $2\pi/(n + 2).$
Tangent at $ p$ to each ray $L_i$, there is exactly one maximal principal line $S_i$ of \fp 2 which
approaches $p$. Two consecutive lines  $S_i$, $ S_{i+1}$, $i = 0,1,2,\ldots ,n + 1 $ $(S_{n+2}  = S_0)$,
bound a hyperbolic sector of \fp 2.
The angular sectors bounded by $ L_i$ and $L_{i+1}$ are bisected by rays $l_i$ , $ i = 0 , 1 , \ldots , n
+ 1$, which play for \fp 1 the same role as $ L_i$ for \fp 2. See Fig.  \ref{fig:umbh} for an illustration. The
lines $ S_i$, are called separatrices of \fp 2  at $ p.$ Similarly, for \fp 1.

\subsection{   Curvature Lines around  Umbilic Curves} \label{sc:uc}
\hskip 2cm

\medskip

In this subsection the results of \cite{gs_cara} will be outlined.

The interest on the structure of principal lines in a neighborhood
of a continuum of umbilic points, forming a curve,  in an analytic
surface goes back to the lecture of Ca\-ra\-th\'eo\-dory \cite{car}.

Let $c:[0,l]\to \mathbb R^3$ be a regular curve parametrized by
arc length $u$ contained in a regular smooth surface $\mathbb M$,
which is
 oriented by the once for all given  positive unit normal vector field
$N$.

Let $T\circ  c =c^\prime$. According  to  Spivak \cite{Sp},  the
Darboux frame $\{T, N\wedge T,N\}$ along  $c$ satisfies the
following system of differential equations:

\begin{equation}\label{eq:da}
\aligned T^\prime &= k_g N\wedge T +k_n N\\
(N\wedge T)^\prime &=-k_g T +\tau_g N\\
N^\prime  &=-k_n T - \tau_g (N\wedge T)\endaligned
\end{equation}

\noindent where $k_n$ is the {\it normal curvature}, $k_g$ is the
{\it geodesic curvature} and $\tau_g$ is the {\it geodesic
torsion} of the  curve $c$.

\begin{proposition} \label{prop:canaleta} Let $c: [0,l]\rightarrow
\mathbb M$ be a regular arc length parametrization of a curve of
umbilic points,  such that $\{ T,N\wedge T, N\}$ is a positive
frame of ${\mathbb R}^3$. Then the expression
{\small \begin{equation}\label{eq:canaleta} \alpha (u,v)= c(u) + v   (N\wedge T)(u)
+ [\frac 12
  k (u)v^2 +\frac 16 a(u) v^3+ \frac 1{24} b(u)v^4+ h.o.t]N(u),
\end{equation}
}
\noindent where  $k(u)=k_n(c(u),T)=k_n(c(u),N\wedge T) $ is the
normal curvature of $\mathbb M$ in the directions $T$ and $N\wedge
T$, defines a local $C^\infty $ chart
 in a small tubular
neighborhood of $c$. Moreover $\tau_g(u)=0$.
\end{proposition}

 The differential equation of curvature lines
in the chart $\alpha $  is given by
{\small
 \begin{equation}\label{eq:cl}
  \aligned (F&g-Gf)dv^2+(Eg-Ge)dudv+(Ef-Fe)du^2=\\
  =&Ldv^2+Mdvdu+Ndu^2=0;
   \endaligned\end{equation}
  $$\aligned  L=&-[ k^\prime v+\frac 12(k_g k^\prime +a^\prime)
 v^2+ \frac 16(k_g a^\prime+ 3k^\prime k_g^2+ b^\prime+
3k^2k^\prime)v^3+ O(v^4) ]  \\
 M=&\;\;  a(u)v+\frac 12[b(u)-3k^3-k^{\prime\prime}-3k_g a(u)]v^2\\
+&  \frac 16[ 15k^3k_g - 3k_g^\prime k^\prime +
(3k_g^2-16k^2)a(u)-a^{\prime\prime} -5k_g b(u) ]v^3+O(v^4) \\
 N=& \;k^\prime v+\frac 12(a^\prime-3k_g k^\prime)
 v^2+ \frac 16(3k^\prime k_g^2 -9k^2k^\prime -5k_g
a^\prime+b^\prime)v^3+
 O(v^4)
 \endaligned$$
 }

\begin{figure}[htbp]
\begin{center}
\includegraphics[angle=0, height=2.5cm]{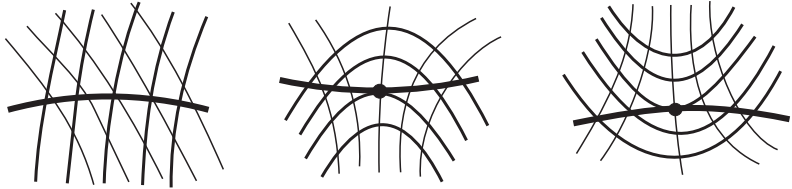}
\end{center}
\caption{ \label{fig:1c} Principal curvature lines near an umbilic
curve: transversal case, left,  and tangential  case, center and
right.}
\end{figure}

\begin{proposition}\label{prop:1c} Suppose that $\nabla {\mathcal
H}(u,0)=(k^\prime, a(u)/2)$ is not zero at a point $u_0$.
 Then the
principal foliations near the point  $c(u_0)$
are as
follows.

\begin{itemize}
\item[i)] If $k^\prime(u_0)\ne 0$ then both principal foliations
are transversal to the curve of umbilic points. See Fig.
\ref{fig:1c}, left.

\item[ii)] If $k^\prime(u_0)= 0$,  $k^{\prime\prime}(u_0)\ne  0$
and $a(u_0)\ne 0$, then one principal foliation is transversal to
$c$ and the other foliation has quadratic contact with the curve
$c$ at the point $c(u_0)$. See Fig.  \ref{fig:1c}, center and right.
\end{itemize}
\end{proposition}

\begin{proposition}\label{prop:11} Suppose that
$k^\prime(0)=a(0)=0 $, $a^\prime(0)k^{\prime\prime}(0)\ne 0$,  at
the point   $c(0)$ of a regular curve $c$ of umbilic points. Let
$A:= -2k^{\prime\prime}(0)/a^\prime (0)\ne 0$ and $B:=
[b(0)-3k(0)^3-k^{\prime\prime}(0)]/a^\prime (0)$. Let $\Delta$ and
$\delta$ be defined by
$$\Delta= -4 A^4+12 B A^3-(36+12 B^2) A^2+(4 B^3+72 B) A-9
B^2-108; \; \; \delta =2- AB.$$

Then the principal foliations at this point are as follows.

\begin{itemize}

\item[i)] If   $ \delta <0 $ and $\Delta <0$ then $0$ is
topologically equivalent to  a Darbouxian umbilic of type  {\rm D}$_1$,
through  which the umbilic curve is adjoined  transversally  to
the
 separatrices.
 See Fig.
\ref{fig:dl} left.

\item[ii)] If   $ \delta  <0 $ and $\Delta>0$ then $0$ is
topologically equivalent to a Darbouxian umbilic of type  {\rm D}$_2$,
through which the  umbilic
 curve is adjoined,  on the interior of the parabolic sectors,
transversally  to the
 separatrices and to  the nodal central line.
See     Fig.  \ref{fig:dl} center.

\item[iii)] If $  \delta >0$ then  $0$ is topologically  a
Darbouxian umbilic of type  {\rm D}$_3$,   through which the umbilic
curve is adjoined transversally  to the
 separatrices.
See   Fig.  \ref{fig:dl} right.
\end{itemize}

\begin{figure}[htb]
\begin{center}
\includegraphics[angle=0, height=3.0cm]{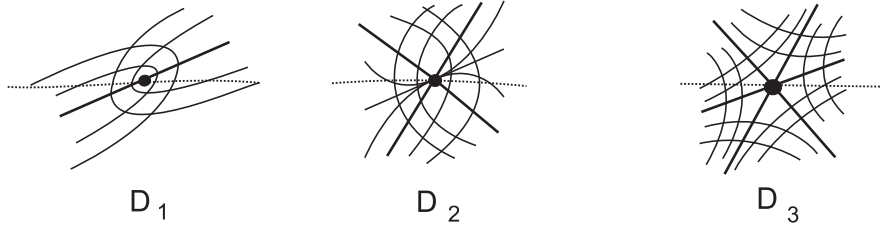}
\end{center}
\caption{ \label{fig:dl} Principal curvature lines near a
Darbouxian-like point on an  umbilic curve, dotted.}
\end{figure}
\end{proposition}

In the previous propositions  have been studied a sample of  the most
generic situations, under the restriction on  the surface  of
having an umbilic curve. Below will be considered the case where
$k$ is a constant which implies the additional constrain that the
umbilic curve be spherical or planar, a case also partially
considered in  Carath\'eodory \cite{car}. Under this double imposition
the simplest patterns of principal curvature lines are analyzed in
what follows.

\begin{figure}[hb]
\begin{center}
\includegraphics[angle=0, height=2.5cm]{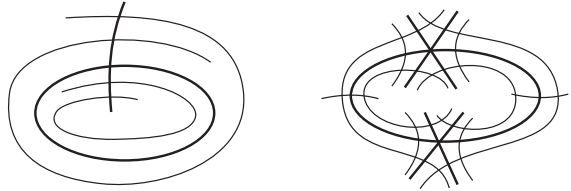}
\end{center}
\caption{ \label{fig:2c} Curvature lines near a spherical umbilic
curve.  }
\end{figure}

\begin{proposition}\label{prop:4c} Let $c$ be a regular closed spherical
 or planar
curve. Suppose that $c$ is a regular curve of umbilic points on a
smooth surface.  Then  the principal foliations near the curve are
as follows.

\begin{itemize}
\item[i)] If ${\mathcal H}_v(u,0)=a(u)/2 \ne 0$ and  $a(u)>0$ for
definiteness,  then one principal foliation is transversal to the
curve $c$ of umbilic points.

The other foliation defines a first return map (holonomy) $\pi$
along the oriented umbilic curve $c$, with first derivative
$\pi^{\prime}=1$ and second derivative given by a positive
multiple of
 $$\int_0^l [k_g(u){a^\prime(u)}/{a(u)^{\frac 32}}]du.$$
  When the above integral is non zero the principal lines  spiral
towards or away from $c$,
   depending on their side relative to
$c$.

\item[ii)] If $a(u)$ has only transversal zeros, near them    the
principal foliations have the topological behavior of a Darbouxian
umbilic point {\rm D}$_3$ at
 which a separatrix has been replaced  with the umbilic curve.
See Fig.  \ref{fig:2c}.
\end{itemize}
\end{proposition}

\section{  Curvature Lines in the Neighborhood of  Critical Points} \label{sc:sp}

In this section, following
\cite{ggsto}, \cite{gscone} and \cite{gstr},  will be described the local behavior of principal curvature  lines  near critical points of the surface such as  Whitney umbrella critical points, conic  critical points and elementary ends of immersions  with constant mean curvature.

\subsection{ Curvature Lines around Whitney Umbrella Critical Points} \label{sc:wsp}
\hskip 2cm

\medskip

In this subsection it will be described the behavior of principal lines near critical points of Whitney type of smooth immersions $\alpha:\mathbb M\to \mathbb R^3$.

The mapping $\alpha$ is said to have a Whitney umbrella at $0$
provided it has rank 1 and its first jet extension $j^1\alpha$ is transversal
to the codimension 2 submanifold $S^1(2,3)$ of 1-jets of rank 1 in
the space $J^1(2,3)$ of 1-jets
of smooth mappings of $(\mathbb R^2, 0)$
to $(\mathbb R^3,0)$.
In coordinates this means that there exist a local chart $(u,v)$
such that $\alpha_u(0)\ne 0,\;\; \alpha_v(0)=0$ and
$[\alpha_u,\alpha_{uv},\alpha_{vv}]\ne 0$. Here $[.,.,.]$ means the determinant of
three vectors.

\begin{figure}[htbp]
\begin{center}
\includegraphics[angle=0, height=3.0cm]{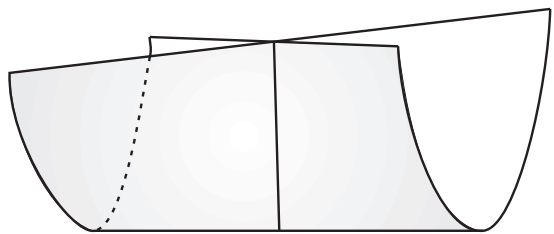}
\end{center}
\caption{ \label{fig:ww} Whitney umbrella critical point.}
\end{figure}

 The structure of a smooth map near such point is illustrated in
Fig.  \ref{fig:ww}. It follows from
the work of Whitney \cite{Wh1} that these points are isolated
and in fact have the following  normal form under diffeomorphic changes
of coordinates in the source and target ($\mathcal A$-equivalence):
 $$x=u,\quad y=uv,\quad z= v^2.$$

 For the study of principal lines the following proposition is useful.

\begin{proposition} \label{prop:ww} Let $\alpha:(\mathbb R^2,0)\to (\mathbb R^3,0)$ be
 a $C^r,\; r\geq 4,$ map with a  Whitney umbrella at $0$.
Then by the action of the group
 ${\mathcal G}^k$ and that   of   rotations  and   homoteties
of $\mathbb R^3$, the map $\alpha$ can be written in  the following form:
$$
\alpha(u,v)=(u,y(u,v),z(u,v))
$$
where,
{\small $$
\aligned
y(u,v) &= uv+ \frac a6  v^3+O(4) \\
z(u,v) &=  \frac b2 u^2+cuv+v^2+
        \frac A6  u^3
+ \frac B2 u^2v+ \frac C2  uv^2+  \frac D6  v^3+O(4)
\endaligned
$$
}
\noindent and $O(4)$ means terms of order greater than  or equal to four.
\end{proposition}

\begin{figure}[htbp]
\begin{center}
\includegraphics[angle=0, height= 4.0cm ]{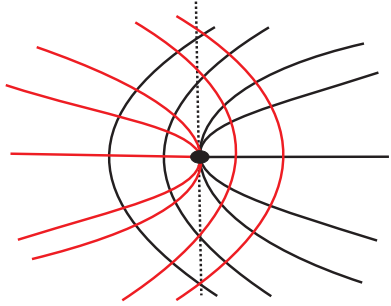}
\end{center}
\caption{ \label{fig:wh}   Curvature lines near a Whitney umbrella.
  }
\end{figure}

 The differential equation of the lines of
curvature of the map $\alpha$ around a Whitney umbrella
 point $(0,0)$,
as in Proposition \ref{prop:ww}, is given by:
{\small
$$\aligned
& [8 v^3 + u \;O(\sqrt{u^2 +v^2}) + O(u^2 + v^2)]dv^2 \\
& +[2u +Cu^2+(D-ac)uv-av^2+ O((u^2 + v^2)^{3/2})]dudv \\
& +[-2v+\frac B2u^2+\frac 12(ac-D)v^2 - b^2 c u^3 + \\
& \qquad \qquad \qquad v \; O(\sqrt{u^2 +v^2}) + O(u^2 + v^2)]du^2=0.
\endaligned
$$
}
\begin{theorem}\cite{ggsto} \label{th:sw} Let $p$ be a Whitney umbrella
 of a map
$\alpha:\mathbb M^2\to \mathbb R^3$ of class $C^k,\;\; k\geq 4$.
Then the principal configuration near
$p$  has the following structure:
Each principal foliation ${\mathcal F}_{i}(\alpha)$ of $\alpha$ has
exactly two sectors at $p;$
one parabolic  and the other
hyperbolic. Also, the separatrices
 of these sectors are tangent to the kernel of $D\alpha_p$.

Fig.  \ref{fig:wh}  illustrates the behavior of
principal curvature lines near a  Whitney umbrella.
 \end{theorem}

\begin{remark} Global aspects of principal configurations  of maps with Whitney umbrellas was carried out in
\cite{ggsto}.
\end{remark}

\medskip
\subsection{  Curvature Lines near  Conic Critical Points} \label{sc:prlc}
 \hskip 2cm

\medskip
A surface in Euclidean $\mathbb R^3$-space is defined implicitly as the variety ${\mathcal V}(f)$ of zeroes of a real valued function $f$, assumed of class $C^k$.

The points $p\in {\mathcal V}(f)$ at which the first derivative $df_p$ does not vanish (resp. vanishes) are called {\it regular} (resp. {\it critical}); they determine the set denoted ${\mathcal R}(f)$ (resp.  ${\mathcal C}(f)$), called the {\it regular} or {\it smooth} (resp. {\it critical}) part of the surface.

The {\it orientation} on  ${\mathcal V}(f)$, or rather
${\mathcal R}(f)$,
 is defined by taking the gradient $\nabla f$ to be the {positive normal.} In canonical coordinates $(x,y,z)$ it follows that $\nabla f=(f_x,f_y,f_z)$.

The Gaussian {\it  normal map} $N_f$, of  ${\mathcal R}(f)$ into $\mathbb S^2$ is defined by $N_f=\nabla f/|\nabla f|$.

The eigenvalues $-k^1_f(p)$ and $-k^2_f(p)$ of $DN_f(p)$, restricted to $T_p {\mathcal R}(f)$, the tangent space to the surface at $p$, define the {\it principal curvatures}, $ k^1_f(p)$ and $ k^2_f(p)$ of the surface at the point $p$. It will be assumed that $k^1_f(p)\leq k^2_f(p)$.

The points on ${\mathcal R}(f)$ where the two principal curvatures coincide, define the set ${\mathcal U}(f)$ of umbilic points of ${\mathcal V}(f)$.

On ${\mathcal R}(f)\setminus {\mathcal U}(f)$, the eigenspaces of $DN_f$ associated to $-k^1_f(p)$ and $-k^2_f(p)$ define $C^{k-2}$ line fields   ${\mathcal L}_1(f)$ and ${\mathcal L}_2(f)$, mutually orthogonal, called respectively {\it minimal} and {\it maximal principal lines fields} of the surface  ${\mathcal V}(f)$. Their integral curves are called respectively the {\it lines of minimal} and {\it maximal principal curvature} or, simply the {\it principal lines} of  ${\mathcal V}(f)$.

The integral foliations  ${\mathcal F}_1(f)$ and  ${\mathcal F}_2(f)$ of the lines fields  ${\mathcal L}_1(f)$ and  ${\mathcal L}_2(f)$ are called, respectively, the {\it minimal} and {\it maximal foliations} of  ${\mathcal V} (f)$.

The net ${\mathcal P} (f)=({\mathcal F}_1(f),{\mathcal F}_2(f))$  of orthogonal curves on ${\mathcal R}(f)\setminus {\mathcal U}(f)$ is called {\it principal net}.

A surface  ${\mathcal V} (f)$ of class $C^k,\; k\geq 3$, is said to have a {\it non degenerate critical point}
at $p$ provided that function $f$ vanishes together with its first partial derivatives at $p$ and that the determinant of the Hessian matrix
$$ H(f,p)=\left(\begin{matrix} f_{xx} & f_{xy} & f_{xz} \\
f_{xy} & f_{yy} & f_{yz}\\
f_{xz} & f_{yz} & f_{zz}\end{matrix}\right)$$

\noindent does not vanish.

The local differentiable structure of ${\mathcal V} (f)$, is determined, modulo diffeomorphism, by the number $\nu=\nu(f,p)$ of negative eigenvalues of $H(f,p)$; $\nu$ is called the {\it index} of the critical point.

It will be assumed that $\nu=1$ and that in the orthonormal coordinates $(x,y,z)$ the critical point is $0$   such that the  diagonal quadratic part of $f$ is given by: $f_2(x,y,z)= x^2/a^2+y^2/b^2-z^2$.

So, $f=f_2+f_3$ with
$$f_3(x,y,z)=\sum_{i+j+k=3}^{} a_{ijk}(x,y,z) x^i y^j z^k,$$

\noindent where $a_{ijk}$ are functions of class $C^{k-3}$.

\begin{figure}[htbp]
\begin{center}
\includegraphics[angle=0, height= 4.7cm ]{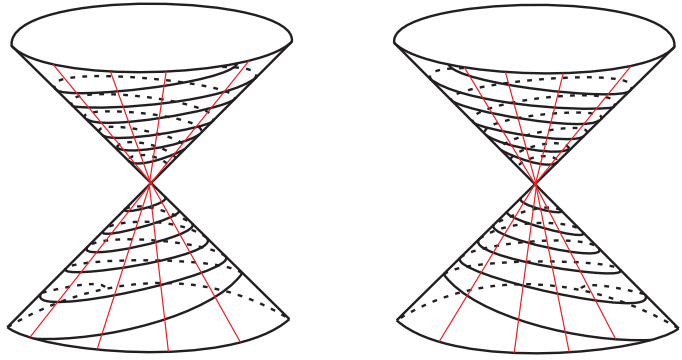}
\end{center}
\caption{ \label{fig:clc} Curvature lines near a critical point of conical type.
  }
\end{figure}

\begin{theorem}\cite{gscone}  Let $f_2(x,y,z)=x^2/a^2+y^2/b^2-z^2$ and
$$  f_{111} (x,y,z)= f_2(x,y,z) + cxyz = 0, \quad \sigma=(a-b)c\ne 0.$$

\noindent Then   ${\mathcal F}_1(f_{111})$ spirals locally around the critical conic point. See Fig.  \ref{fig:clc}.

 The other principal foliation  ${\mathcal F}_2(f_{111})$ preserves locally the radial behavior of  ${\mathcal F}_1(f_{2})$.
 More precisely, there is a local orientation preserving homeomorphism  mapping   $\{{\mathcal V} (f_{111}),0\}$ to
 $\{{\mathcal V} (f_{2}),0\}$, sending the principal net $ {\mathcal P} (f_{111}) $ on the net  ${\mathbb F}=({\mathbb F_1}, {\mathbb F_2 })$ on ${\mathcal V}(f_2)$, defined by the integral foliations ${\mathbb F_1}$ and ${\mathbb F_2}$ of the vector fields
 $$X_1=X_2- a^2 y\partial/{\partial x} +b^2 x \partial/{\partial y}, X_2=\sigma z(x\partial/{\partial x}+y\partial/{\partial y}+z\partial/{\partial z}).$$

\end{theorem}

\begin{theorem} \cite{gscone}  Let $f$ be of class $C^k$, $k\geq 7$, with a non degenerate critical point of index $\nu=1$ at $0$, written as $f=f_2+f_3$. For $\sigma=(a-b)a_{111}(0)\ne 0$, there is a local orientation preserving homeomorphism  mapping   $\{{\mathcal V} (f),0\}$ to  $\{{\mathcal V} (f_{111}),0\}$,
  sending the principal net $ {\mathcal P} (f) $ on that of $ {\mathcal P} (f_{111}). $
  Here $f_{111}(x,y,z)=f_2(x,y,z)+a_{111}(0) xyz$.

\end{theorem}
\medskip

\begin{remark} Local and global stability of principal nets  ${\mathcal P} (f)$ on surfaces defined implicitly were also studied in \cite{gscone}.\end{remark}

\subsection{ Ends of Surfaces Immersed with Constant Mean Curvature} \label{sc:endh}
\hskip 2cm

\medskip

 Let $(u ,v) :\mathbb M \to \mathbb R^2\setminus \{0\} $ be isothermic coordinates for an
immersion $\alpha :\mathbb M \to\mathbb  R^3$ with constant mean curvature. If the associated complex
function $\phi(w),\; w = u + iv, $ has in zero a pole of order $n \ne  2$, then there exists a
small neighborhood $V$ of $0$ in $\mathbb R^2$ such that the principal lines of $\alpha$, restricted to
$(u ,v)^{-1}  (V\setminus\{ 0\}) $, distribute themselves (modulo topological equivalence) as if the
associated complex function were $w^{-n}$.

An {\it end } of the
immersion $\alpha$
defined by the system of open sets $U_j=\{(u,v): u^2+v^2<1/j, \; j\in \mathbb N\}$
where $(u , v)$ are isothermic coordinates for $ \alpha$, on which the associated complex
function $\phi$ has a pole of order $n $ in $(u , v) = (0,0)$, is called an {\em elementary end}  of
order $n$ of  $\alpha$. The index of such an elementary end is $n/2.$

\begin{figure}[htbp]
\begin{center}
\includegraphics[angle=0, height= 6.5cm ]{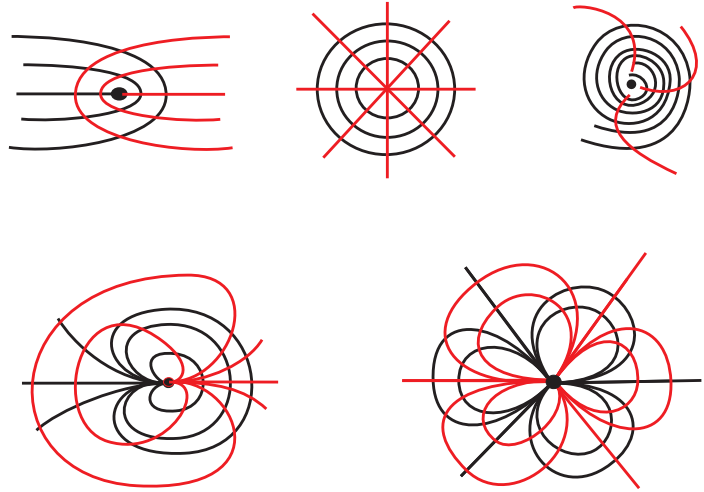}
\end{center}
\caption{ \label{fig:poloh} Curvature lines near elementary end points.
  }
\end{figure}

The lines of curvature of an immersion $ \alpha$   with constant mean
curvature near an elementary end $E$ of order $n$ are described as follows. See more details in \cite{gstr} and Fig.  \ref{fig:poloh}.

\begin{itemize}
\item[a)] For $n = 1$, there is exactly one line $S$ (resp. $s $) of \fp 1 (resp. \fp 2) which tends to
$E$, all the other lines fill a hyperbolic sector bounded by $E$ and $S$ (resp. $s$).

\item[b)] For $ n = 2,$ suppose that $\phi(z)$  is the associated complex and $ a = \lim_{z\to 0} z^2\phi(z).$
There are two cases:
\item[b.1)]  $a \notin \mathbb R \cup (i\mathbb  R)$. Then the lines of \fp 1 and \fp 2, tend to E.
\item[b.2)] $a \in \mathbb R \cup (i\mathbb  R)$.  Then the lines of \fp 1 (resp. \fp 2) are circles or rays tending to $E$
and those of \fp 2, (resp. \fp 1) are rays or circles.
\item[c)] For $n \geq 3$, every line of \fp 1 and \fp 2, tends to $E.$ The principal lines distribute
themselves into $n - 2$ elliptic sectors, two consecutive of which are separated by a
parabolic sector.
\end{itemize}

\section{ Curvature Lines near Principal Cycles}\label{sc:pclpc}

A compact leaf $\gamma$ of \fp 1 (resp. \fp 2) is called  a {\it
minimal} (resp. {\it ma\-xi\-mal}) {\it principal cycle.}

A useful local parametrization near a principal cycle is given by the following proposition and was introduced by Gutierrez and Sotomayor in \cite{gs1}.

\begin{proposition} \label{prop:canal} Let $\gamma: [0,L]\rightarrow {\mathbb R}^3 $ be a principal cycle  of an immersed surface $\mathbb M$ such that
$\{ T,N\wedge T, N\}$ is a positive frame of ${\mathbb R}^3$. Then the expression
\begin{equation}\label{eq:canal} \aligned \alpha (s,v)=& \gamma(s) + v   (N\wedge T)(s) \\
+& [\frac 12
  k_2(s)v^2 + \frac 16 b(s) v^3+
    o(v^3)]N(s),
 \;\; -\delta < v <\delta\endaligned
\end{equation}

\noindent where  $k_2$ is the principal  curvature   in the direction of $N\wedge T$,  defines a local $C^\infty $ chart on the surface $ {\mathbb M}$  \; defined in
a small tubular neighborhood of $\gamma$.
\end{proposition}

 \begin{remark}

 Calculation shows that the following relations hold

 \begin{equation}\label{eq:kgo}
 k_g(s)=- \frac{(k_1)_v}{k_2-k_1}, \;\;k_g^\perp(s)=- \frac{(k_2)^\prime}{k_2-k_1},\;\;\;\;\; b(s)=(k_2)_v=\frac{\partial k_2}{\partial v}
 \end{equation}
\noindent Here $k_g^\perp(s)$ is the geodesic curvature of the maximal principal curvature line which pass through $\gamma(s)$.

 \end{remark}

\begin{proposition}[Gutierrez-Sotomayor]\label{prop:frm}\cite{gs1} Let $\gamma$ be a  minimal principal cycle of an immersion $\alpha:\mathbb M\to \mathbb R^3$ of length $L$. Denote by $\pi_\alpha$ the first return map associated to $\gamma$. Then

\begin{equation}
\aligned
\pi^\prime_\alpha(0)=&exp [ \int_\gamma\frac{-dk_2}{k_2-k_1}]= exp [\int_\gamma k_g^\perp(s)ds]\\
=& exp [\int_\gamma\frac{-dk_1}{k_1-k_2}]=exp [\frac 12\int_\gamma \frac{d{\mathcal H}}{\sqrt{{\mathcal H}^2-{\mathcal K}}}].\endaligned
\end{equation}

\end{proposition}

The following result established in \cite{gs3} is improved in the next  proposition.

 \begin{proposition} \label{prop:sh} Let $\gamma$ be a  minimal principal cycle of length $L$ of a surface $ \mathbb M\subset \mathbb R^3$. Consider a chart $(s,v)$ in a neighborhood of $\gamma$ given by equation \eqref{eq:canal}.  Denote by $k_1$ and $k_2$ the principal curvatures of $\mathbb M$. Let $Jac(k_1,k_2)=\frac{\partial(k_1,k_2)}{\partial(s,v)}=(k_1)_s (k_2)_v - (k_1)_v (k_2)_s $ and suppose that $\gamma$ is not hyperbolic, i.e.   the first derivative of the first return map $\pi$ associated to $\gamma$ is one. Then the  second derivative of $\pi$ is given by:
$$  \pi^{\prime\prime}(0)= \int_0^L e^{ -\int_0^s \frac{k_2^\prime}{k_2-k_1} du} \frac{Jac(k_1,k_2)}{(k_2-k_1)^2 }ds.  $$
\end{proposition}

  \begin{theorem}\label{th:cn}  \cite{gspit,gsagag} For a principal cycle $\gamma$ of multiplicity $ n$,
 $ 1 \le n < \infty$, there is
a chart $ (u , v)$ with $ c = \{ v = 0 \}$,  such that the
principal lines are given by

$$\aligned   du =& 0,\;\; dv - a_1 v du,\;  a_{1}\in \mathbb R, \;\text{for}\; n=1.\\
 du =& 0,\;\;dv - v^n (a_n -{a}_{2n -1} v^{n-1})
du = 0 ;\;\; a_{n} ,  a_{2n-1}\in \mathbb R, \;\text{and }\; n\geq 2.\endaligned$$

Here, the numbers $a_1$, $a_n$ and ${a}_{2n-1}$ are uniquely determined by the
jet of order $2n+1$  of immersion $\alpha$ along $\gamma$ and
can be expressed in terms of   integrals involving the principal
  curvatures $k_1(u,v) $ and $k_2(u,v)$ and its derivatives.

  \end{theorem}

 The next result  establishes how a principal cycle of multiplicity $n$,
$ 2 \le n < \infty,$ of
an immersion $\alpha$ splits under deformations $\alpha_{\epsilon}$.

Recall from \cite{GoSc}   that a family of functions $\;
U(.,\epsilon),\;$ $ \epsilon \in \mathbb R^n,$ is an universal unfolding
of $U(.,0) $ if for any deformation
$H(.,\theta),$   $\theta \in  \mathbb R^l$ of $ U(.,0), $
 the following  equation
 $ \; H(v,\theta) = S(v,\theta)
U(\beta(v,\theta),\lambda(\theta))\; $ holds.

 Here $ \;S,\;$  $\beta\; $ and $\,\lambda\,$ are
 $C^{\infty}$
functions with $\, S(v,\theta) > 0,\,$  $\beta(v,0) = v$,
 $\; \lambda(0)=0$.
Furthermore it is required that $ n $ is the minimal number
 with this property.

For a minimal principal cycle $\gamma$ of $\alpha$ on which $k_1$ is not constant
consider the following deformation:

$${\alpha}_{\epsilon} = \alpha + k_1^{\prime}(u) \delta(v) \big(\sum_{i=1}^{n-1}
{\epsilon}_i \frac{v^i}{i!} \big) N_{\alpha} (u).$$

Here $v $ is a $C^{\infty}$ function on M such that $ v{\vert}_c =
0$, ${\nabla}_{\alpha} v $ is a unit vector on the induced metric ${< ,
>}_{a}$; $ {\nabla}_{\alpha}$ denotes the gradient relative to this metric
and $ \delta $ is a non negative function,
identically $1$ on a neighborhood of $ \gamma, $ whose
support is contained on the domain of $ v$.

\begin{theorem} \label{th:duniv}\cite{gsagag} For a minimal principal cycle $ \gamma $ of multiplicity
n, $ 2 \le n < \infty,$ of $\alpha$, on which the principal curvature $ k_1$ is not a
constant, the following holds.

 The function $ U(x,\epsilon) = {\pi}_{{\alpha}_{\epsilon}} (x) - x $
provides a universal unfolding for  $ U(.,0) = {\pi}_{\alpha} - id $.
Here ${\pi}_{{\alpha}_{\epsilon}}
$ denotes the return map of the deformation ${\alpha}_{\epsilon}$.
 \end{theorem}

\begin{remark}\begin{itemize}

\item[i)] Principal cycles on immersed surfaces with constant mean
curvature was studied in \cite{gstr}. There is proved that the
Poincar\'e transition map preserves a transversal measure   and the principal
cycles  appears in open sets, i.e.  they fill an open region.

\item[ii)] Principal cycles on immersed   Weingarten surfaces have
been  studied in \cite{soagag}.  They also appear in open sets.

    \item[iii)] In
\cite{gspit} was established an integral expression for  $\pi^{\prime}$ in
terms of the principal curvatures and the  Riemann Curvature Tensor of
the manifold in which the surface is immersed. It seems challenging to
discover the general pattern for the higher derivatives of the return
map in this case.
\end{itemize}
\end{remark}

\section{Curvature  Lines on Canal Surfaces} \label{sc:cs}

In this section it will   determined the principal curvatures and
principal curvature lines on {\it canal surfaces} which are the
envelopes of families of spheres with variable radius and centers
moving along  a closed regular curve in $\mathbb R^3$, see \cite{ve}.
This study
were carried out in \cite{glsabc}.
\medskip

Consider the space $\mathbb R^3$  endowed with the Euclidean inner
product $<\, , \,>$ and  norm $|\,\,\,| = <\, ,\,>^{1/2}$ as well
as with a canonical orientation.

Let ${\mathbf c}$ be a smooth regular closed curve immersed  in
$\mathbb R^3$, parametrized by arc length $s\in [0,L]$. This means
that
\begin{equation} \label{eq:ct}
{\mathbf c}^{\prime} (s)={\mathbf t}(s), \; |{\mathbf t}(s)|=1,
\;{\mathbf c}(L)={\mathbf c}(0)\;  .
\end{equation}
Assume also that the curve is {\it bi-regular}. That is:
\begin{equation} \label{eq:bireg}
\kappa(s) = |{\mathbf t}^{\prime} (s)| >0.
\end{equation}

Along ${\mathbf c}$ is defined  its  moving {\it Frenet frame}
$\{{\mathbf t},{\mathbf n},{\mathbf b}\}$. Following  Spivak \cite{Sp}
and Struik \cite{St}, this  frame is {\it positive}, {\it
orthonormal} and verifies {\it Frenet equations}:
\begin{equation} \label{eq:fr}
{\mathbf t}^{\prime} (s) = \kappa(s){\mathbf n}(s), \; {\mathbf
n}^{\prime} (s) = -\kappa(s){\mathbf t}(s) + \tau (s) {\mathbf
b}(s), \; {\mathbf b}^{\prime} (s) = -\tau (s){\mathbf n}(s).
\end{equation}

Equations (\ref{eq:ct}) to  (\ref{eq:fr}) define the {\it unit
tangent}, ${\mathbf t}$, {\it principal normal}, ${\mathbf n}$,
{\it curvature},  $\kappa$,  {\it binormal}, ${\mathbf b}={\mathbf
t}\wedge {\mathbf n}$, and {\it torsion},  $\tau$,  of the
immersed curve ${\mathbf c}$.

\begin{proposition}\label{prop:immer}
Let $r(s)>0$ and $\theta(s) \in \; ]0,\pi[$ be smooth functions of
 period  $L$. The mapping $\alpha:\mathbb T^2 = \mathbb S^1\times
\mathbb S^1 \to\mathbb R^3 ,$ defined on $\mathbb R^2$ modulo
$L\times 2\pi$ by
\begin{equation}\label{eq:al}
\alpha(s, \varphi)={\mathbf c}(s)+ r(s)\cos\theta(s) {\mathbf
t}(s)+  r(s)\sin\theta(s)[\cos \varphi\, {\mathbf n}(s)+
\sin\varphi\, {\mathbf b}(s)],
\end{equation}
is tangent to the sphere of center ${\mathbf c}(s)$ and radius
$r(s)$ if and only if
\begin{equation}\label{eq:r1}
\cos\theta(s) = -r^{\prime} (s).
\end{equation}
Assuming \eqref{eq:r1}, with $r^{\prime} (s) < 1$, $\alpha$ is an
immersion provided
\begin{equation}\label{eq:reg1}
\kappa(s) <\frac{1-r'(s)^2 - r(s)r''(s)}{r(s)\sqrt{1-r'(s)^2}}.
\end{equation}
\end{proposition}

\begin{definition} \label{def:env}
A mapping such as $\alpha$, of $\,\, {\mathbb T}^2$ into $\mathbb R^3$, satisfying
conditions  \eqref{eq:r1} and \eqref{eq:reg1}  will be called an
{\it immersed canal surface} with {\it center along} ${\mathbf
c}(s)$ and {\it radial function} $r(s)$. When $r$ is constant, it
is called an {\it immersed tube}. Due to the tangency condition
\eqref{eq:r1}, the  {\it immersed canal surface} $\alpha$ is the
{\it envelope} of the family of spheres of radius $r(s)$ whose
centers range along the curve  ${\mathbf c}(s)$.
\end{definition}

\begin{theorem}\label{th:env}
Let $\alpha:{\mathbb S}^1\times {\mathbb S}^1 \to{\mathbb R}^3$ be a
smooth immersion expressed  by \eqref{eq:al}. Assume the
regularity conditions \eqref{eq:r1} and \eqref{eq:reg1} as in
Proposition \ref{prop:immer} and also that

\begin{equation} \label{eq:nu}
k(s)<\Big| \frac{1-r'(s)^2-2r(s)r^{\prime\prime}(s) }{2 r(s)
\sqrt{ 1-r'(s)^2}}\Big|.
\end{equation}
The maximal principal curvature lines are the circles tangent to
$\partial/\partial \varphi$. The maximal principal curvature is
\[
k_2 (s) = 1/r(s).
\]
The minimal principal curvature lines are the curves tangent to
\begin{equation}\label{eq:clc}
V(s,\varphi)= \frac{\partial}{\partial s} - \left(\tau(s)+
\frac{r'(s)}{(1-r'(s)^2)^{1/2}}\kappa(s)\sin \varphi\right)
\frac{\partial}{\partial \varphi}.
\end{equation}
The expression
\begin{equation}\label{eq:d}
r(s)(1-r'(s)^2)^{1/2}\kappa(s)\cos\varphi
+r(s)r^{\prime\prime}(s)- (1-r'(s)^2)
\end{equation}
is negative, and the minimal  principal curvature is given by

\begin{equation}\label{eq:k1}
k_1 (s,\varphi) = \frac{\kappa(s)(1-r'(s)^2)^{1/2}\cos\varphi
+r^{\prime\prime}(s)}{r(s)(1-r'(s)^2)^{1/2}\kappa(s)\cos\varphi
+r(s)r^{\prime\prime}(s) - (1-r'(s)^2)}.
\end{equation}

There are  no umbilic points for $\alpha$: $k_1 (s,\varphi) < k_2
(s). $

\end{theorem}

\begin{remark}
A consequence of the
Riccati structure for principal curvature lines on canal immersed
surfaces established in Theorem \ref{th:env}  implies that the maximal number  of isolated periodic
principal lines is $2$. Examples of canal surfaces with two
(simple i.e. hyperbolic) and one (double i.e. semi--stable)
principal periodic lines  were developed in \cite{glsabc}.
\end{remark}

\section{ Curvature Lines near Umbilic Connections and Loops }\label{sec:ucloop}

 A principal line  $\gamma$  which is an {\it umbilic  separatrix} of two
 different umbilic
points  $p$, $q$  of  $\alpha$ or twice a separatrix of the same
umbilic point $ p$  of $\alpha$ is called an {\it umbilic
separatrix connection} of  $\alpha$;  in the second case $\gamma$
is also called an {\it umbilic separatrix loop.}
The simplest
bifurcations of umbilic connections, including umbilic loops, as well
as the consequent
appearance of principal cycles will be outlined below, following
\cite{gsd1d1} and \cite{sell}.

\begin{figure}[htp]
 \begin{center}
 \includegraphics[angle=0,height=8cm, width=11cm]{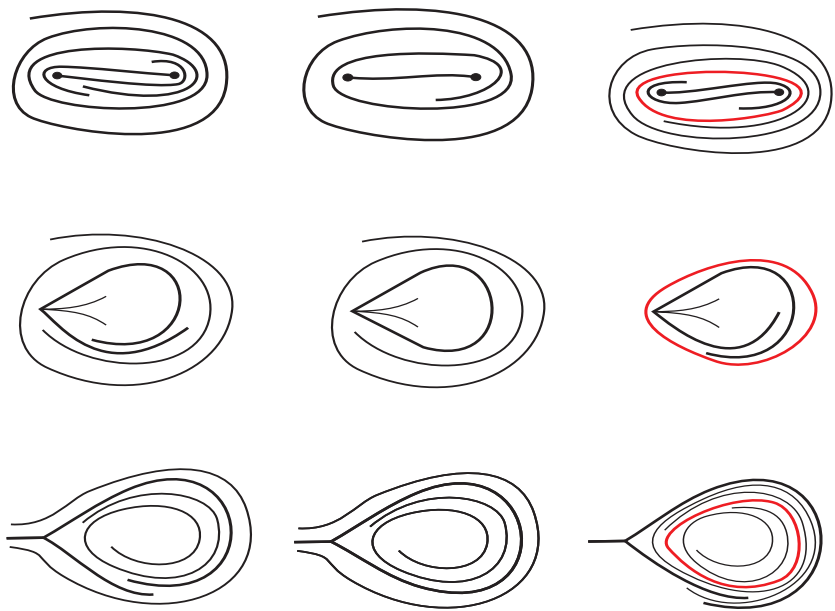}
 \caption{ Simple Connection and Loops and their Bifurcations:   $C_{11}$, top, $C_{22}$,  middle and  $C_{33}$,  bottom. \label{fig:uc}}
   \end{center}
 \end{figure}

There are three types of umbilic connections, illustrated in Fig.  \ref{fig:uc},
from which principal cycles bifurcate.
They are defined as follows:

   $\bullet$  $C_{11}$-{\em simple connection}, which consists in two  {\rm D}$_1$ umbilics joined by their separatrix, whose return map $T$
has first derivative $T^\prime \neq 1$. In \cite{gsd1d1} this derivative is expressed in terms of the third order jet
of the surface at the umbilics\footnote{A misprint in the expression for the asymmetry  $\chi$ at  a  {\rm D}$_1$ point (which, multiplied by $\pi/2$, gives the logarithm of the derivative of the return map at the point)  should be corrected to $\chi =\frac{c}{\sqrt{(a-2b)b-(c^2 /4)}}$. }.
\vspace{.1cm}

    $\bullet$  $C_{22}$-{\em simple loop}, which consists in one {\rm D}$_2$ umbilic point self connected  by a separatrix, whose return map $T$
verifies  $\lim_{x\to 0+}\log(T(x))/\log(x) \neq 1$. In \cite{gsd1d1} this derivative is expressed in terms of the third order jet
of the surface at the umbilic.
\vspace{.1cm}

    $\bullet$  $C_{33}$-{\em simple loop}, the same as above exchanging  {\rm D}$_2$ by  {\rm D}$_3$ umbilic point.

Below   will be outlined the results obtained in \cite{sell}.

There are two bifurcation patterns producing  principal cycles which
 are associated with the bifurcations
of  {\rm D}$^{1}_{2}$  and   {\rm D}$^{1}_{2,3}$  umbilic points, when their
separatrices form loops, self connecting these points. They are defined as
follows.

 \vspace{.1cm}

\begin{figure}[htbp]
 \begin{center} 
 \includegraphics[angle=0, height=3.5cm]{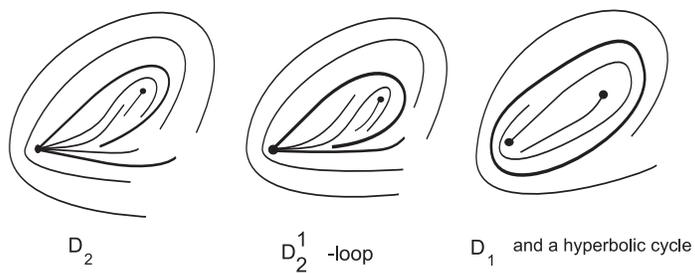}
 \caption{ {\rm D}$^{1}_{2}$ -   loop bifurcation. \label{fig:d12loop}}
   \end{center}
 \end{figure}

 $\bullet$  A  {\rm D}$^{1}_{2}$ - {\it interior loop} consists on a point of type   {\rm D}$^{1}_{2}$ and its
isolated separatrix, which is assumed to be  contained in the interior of
the parabolic sector. See  Fig.  \ref{fig:d12loop}, where such loop together with the
bifurcating principal cycle are illustrated. Here, a hyperbolic principal
cycle bifurcates from the loop when the   {\rm D}$^{1}_{2}$ point bifurcates into   {\rm D}$_1$ .

\begin{figure}[htbp]
 \begin{center}
 \includegraphics[angle=0, height=3.5cm]{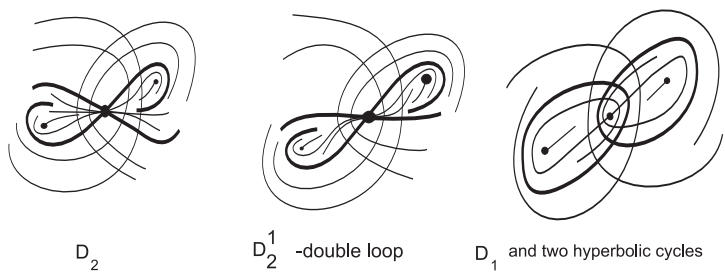}
 \caption{  {\rm D}$^{1}_{2}$ -  double  loop bifurcation. \label{fig:d12loopd}}
   \end{center}
 \end{figure}

    If both principal foliations have    {\rm D}$^{1}_{2}$ - interior loops (at the same
  {\rm D}$^{1}_{2}$   point), after bifurcation there appear two hyperbolic cycles, one for
each foliation. This case will be called {\it double}  {\rm D}$^{1}_{2}$
 - {\it interior loop}. In
Fig.  \ref{fig:d12loopd},
Fig.  \ref{fig:d12loop}
has been  modified and completed accordingly so as to
represent both maximal and minimal foliations, each  with its respective
  {\rm D}$^{1}_{2}$ -interior loops (left) and bifurcating  hyperbolic principal
  cycles (right).

\begin{figure}[htbp]
 \begin{center}
 \includegraphics[angle=0, height=3.0cm]{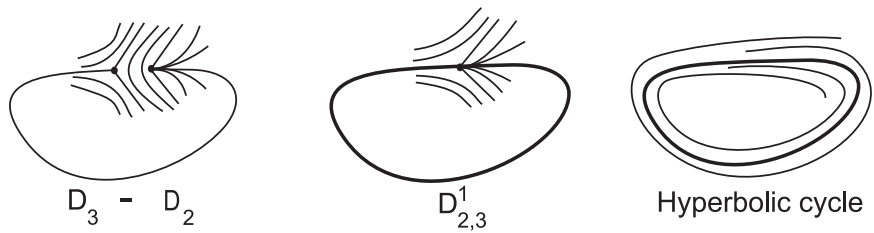}
 \caption{ {\rm D}$^{1}_{2,3}$ - loop bifurcation. \label{fig:d23loop}}
   \end{center}
 \end{figure}

 $\bullet$   A  {\rm D}$^{1}_{2,3}$ - {\it interior loop} consists on a point of type  {\rm D}$^{1}_{2,3}$
     and its
hyperbolic separatrix, which is assumed to be  contained in the interior of
the parabolic sector.
See  Fig.  \ref{fig:d23loop}, where such loop together with the
bifurcating principal cycle are illustrated. Here, a unique hyperbolic principal
cycle bifurcates from the loop when the umbilic points are annihilated.

\section{ Principal Configurations on Algebraic Surfaces in $\mathbb R^3$}\label{sec:as}

\medskip
An  algebraic  surface of degree $n $ in Euclidean $ (x_1, x_2, x_3)$-space $\mathbb R^3$  is defined by the  variety  $A(\alpha)$ of
real zeros of a polynomial $\alpha $ of the form $\alpha= \sum \alpha_h , \,h=0,\,1,\,2,...,n$,
where $\alpha_h$   is  a   homogeneous   polynomial  of  degree  $h$: $
\alpha_h = \sum a_{ijk}x_1^i x_2^j x_3^k, \,i+j+k=h$,  with real
coefficients.

An  {\it end  point}  or {\it point  at infinity} of $A(\alpha) $ is a point
in  the unit sphere $\mathbb S^2$, which is the limit  of  a
sequence  of  the  form $p_n /|p_n|$, for $ p_n$  tending to
infinity in $A(\alpha). $

The  {\it end  locus}  or {\it curve at infinity} of $A(\alpha)$ is defined as
the collection $E(\alpha_n) $ of end  points  of $A(\alpha).$
Clearly, $E(\alpha) $ is contained in  the  algebraic  set $E
(\alpha) = \{p \in \mathbb S^2 ; \alpha (p)= 0\}, \;$ called  the
{\it algebraic  end  locus}  of $A(\alpha).$

A surface $A(\alpha)$ is said to be regular (or smooth) at
infinity if $0 $ is a regular value of the restriction  of $
\alpha_n $  to $\mathbb S^2.$  This is equivalent to require  that
$\nabla \alpha_n (p)\wedge p$ does not  vanish  whenever $\alpha_n
(p)=0,$ $p\ne 0$. In this case, clearly $ E(\alpha) =  E_n
(\alpha)  $ and,  when non empty, it consists of a finite
collection  $\{\gamma_i;\,\,  i=  1,\ldots , k(\alpha_n)\}$  of smooth
closed curves, called the (regular) {\it curves at infinity}  of $
A(\alpha)$; this collection of curves is invariant under  the
antipodal map, $a:p \to -p,$ of the sphere.

There  are  two  types of curves at  infinity:  odd  curves  if $
a(\gamma_i)=\gamma_i$, and even curves if $\gamma_j =a(\gamma_i)$
is disjoint from $\gamma_i$ .

A  regular  end  point  $p$ of $ A(\alpha)$  in $\mathbb S^2$ will
be  called  an ordinary or biregular  end point of $ A(\alpha)$ if
the geodesic curvature, $ k_g$, of  the  curve $E(\alpha) $  at
$p$, considered  as  a  spherical  curve,  is different  from
zero; it is called {\it singular} or {\it inflexion end point} if $k_g$  is
equal to zero.

 An   inflexion   end   point   $p$  of  a  surface  $A(\alpha)$  is
  called {\it bitransversal},
provided the following two transversality  conditions hold:

$$ (T_1) \;\;\;\;\;\;\;\;\;\;\;\;\;\;\;\; k_g^\prime(p) = dk_g (p;\tau) \ne  0$$
and
$$ (T_2) \;\;\;\;\;\;\;\;\; \nu(p) = d[\alpha_{n-1}   /|\nabla \alpha_n |](p;\tau) \ne 0.$$

 There are two
different types of bitransversal inflexion   end  points,
illustrated   in   Fig.  \ref{fig:2}: {\em hyperbolic }   if  $
k_g^\prime .\nu<0  $  and {\em elliptic} if $k_g^\prime .\nu
>0.$

 \begin{figure}[htbp]
 \begin{center}
 \hskip .5cm
  \includegraphics[angle=0, height=2.8cm]{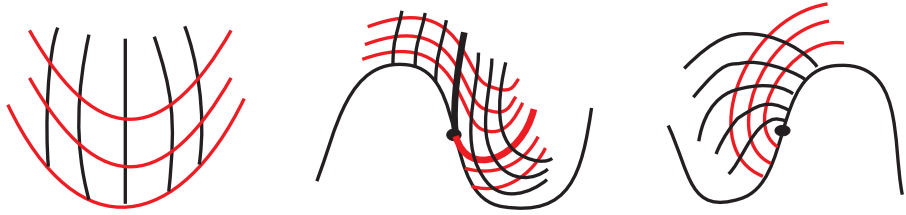}
  \caption{Principal nets near biregular (left) and bitransversal
   end points (hyperbolic, center,  and elliptic, right).
  \label{fig:1a} }
 \end{center}
  \end{figure}

 A regular  end  curve  $\gamma$  all whose  points  are biregular,
i.e. $\gamma$ is free from inflexion points,  will be called a
{\it principal cycle at infinity}; it will be called {\em
semi-hyperbolic} if
 $$ \eta = \int_\gamma \alpha_{n-1}|\nabla\alpha_n|^{-1}d(k_g^{-1}) \ne  0.$$

Little is know  about  the  structure  of  principal  nets  on
algebraic surfaces of degree  $n$.  The  case   of   quadrics  $(n
=  2),$ where  the  principal  nets  are  fully known,   is   a
remarkable  exception, whose study goes back to the classical
works  of Monge  \cite{mo},   Dupin \cite{du}, and
Darboux \cite{da},  among   others.
See the discussion in section \ref{sec:hl} and  also  \cite{Sp}, \cite{St}.

The description of their
principal nets  in terms of  intersections  of the surface with
other families of quadrics  in   triply   orthogonal  ellipsoidal
coordinate   systems,  appear  in  most differential  geometric
presentations of surface theory,  notably  Struik's  and Spivak's
\cite{St}, \cite{Sp}. Also,  Geometry  books  of  general
expository character, such as Fischer \cite{Fi} and   Hilbert-Cohn
Vossen \cite{Hi-CV},  also explain these  properties  of quadrics
and include  pictorial illustrations of their principal nets.

For   algebraic   surfaces   of   degrees   three   (cubics), four
(quartics) and higher,   however,  nothing  concerning  principal
nets, specific  to  their  algebraic character, seems to be known.

Consider the vector space ${\mathcal A}_n $ of all polynomials of
degree  less than or equal to  $n$, endowed with the structure of
${\mathbb R} ^N$ -space
\noindent
defined by the $N=N(n)$ coefficients  of  the  polynomials  $\alpha. $
 Here $ N(n) = (n+1)(n+2)(n+3)/6.$
The
distance in  the  space  ${\mathcal A}_n $  will be denoted $ d_n
(.,.).$

Structural  Stability  for  (the principal net  of)  an  algebraic
surface $A(\alpha) $ of degree $ n$  means  that  there  is an
$\epsilon >0$   such  that for any $\beta$ with $d
(\alpha,\beta)<\epsilon$ there is a homeomorphism  $h$  from
$A(\alpha) $ onto $ A(\beta)$ mapping $U(\alpha)$ to $U(\beta)$
and also mapping the lines of \fp 1  and  \fp 2
onto those of ${\mathcal F}_1(\beta)$ and ${\mathcal F}_2 (\beta)$, respectively.

Denote  by  $\Sigma_n$   the class of surfaces $ A(\alpha), \;
\alpha\in {\mathcal A}_n$, which are regular and regular at
infinity and satisfy simultaneously that:

\begin{itemize}
\item[a)] All its umbilic points are Darbouxian and all inflexion
ends are bitransversal.

\item[b)]  All  its    principal cycles are hyperbolic and
all biregular end curves, i.e. cycles at infinity,  are
semi-hyperbolic.

\item[c)]  There  are  no  separatrix  connections  (outside   the
end locus)  of umbilic and  inflexion end points.

\item[d)]  The  limit  set  of  any principal line is a principal
cycle (finite or infinite), an umbilic point or an end point.
\end{itemize}

  These  conditions extend to algebraic surfaces of degree $n$, $n\geq 3$, the conditions given
 by  Gutierrez  and  Sotomayor  in  \cite{gs1, gsln, gs2}, which imply principal stability for compact
 surfaces.

\medskip

\begin{theorem}\label{th:d}\cite{gsalg}   Suppose that $n\geq 3$.  The  set  $\Sigma_n$   is  open  in
${\mathcal A}_n$, and any surface $A(\alpha) $ on
  it is principally structurally stable.\end{theorem}
\vskip .3cm
\begin{remark}
\begin{itemize}
\item[i)]
 For  $n=2$  the  stable  surfaces  are  characterized,  after  Dupin's  Theorem, by the
 ellipsoids  and hyperboloid of two sheets with different axes and by hyperboloid of one
  sheet (no conditions on the axes). See Theorem \ref{th:q2}.

 \item[ii)] By  approximating,  in  the  $C^3$    topology,   the  compact  principal structurally
 stable  surfaces  of Gutierrez and Sotomayor \cite{gs1, gsln, gs2}  by algebraic ones, can be obtained
 examples   of  algebraic  surfaces  (of  undefined  degree),  which  are  principally
 structurally stable on a compact connected component.

 In this form all the patterns of stable principal configurations of compact smooth
 surfaces are realized by algebraic ones, whose degrees, however,  are not determined.
 \end{itemize}
 \end{remark}
To close this section an open problem is proposed.
\vskip .3cm

\begin{problem} \label{pro:2}

 Determine the class of principally  structurally stable  cubic
and higher degree surfaces. In other words, prove or disprove the
converse of \ref{th:d}.

Prove the density of $\Sigma_n$ in ${\mathcal A}_n$, for $n = 3$
and higher.
\end{problem}

Theorem \ref{th:q2} in section   \ref{sec:2} deals with the case of degree 2
---quadric--- surfaces.
\vskip .3cm

\section{  Axial Configurations on Surfaces Immersed in ${\mathbb R}^4$ } \label{sc:ax}

\medskip
Landmarks of the Curvature Theory for surfaces in ${\mathbb R}^4$
  are the
works of Wong \cite{W} and Little \cite{Lit}, where a review of properties of the Second Fundamental Form,
 the Ellipse of Curvature (defined as the image of this form on unit tangent circles) and related geometric
  and singular theoretic notions are presented.   These authors give a list
of pertinent references  to original sources previous to 1969, to which one must add that of
Forsyth \cite{Fo}.
  Further  geometric properties of surfaces in $\mathbb R^4$ have been pursued by Asperti \cite{As} and
   Fomenko \cite{Fm}, among others.

For an immersion $\alpha$  of a surface  $\mathbb M $ into
${\mathbb R}^4$, the  axiumbilic singularities $ {\mathbb
U}_{\alpha} $, at which the ellipse of curvature degenerates into
a circle,  and the lines of axial curvature are assembled into
two {\it axial configurations}: the {\it principal axial
configuration}: ${\mathbb P}_{\alpha}= \{{\mathbb U}{_\alpha},\;
{\mathbb X}_{\alpha}\}$ and the  {\it mean  axial configuration}:
${\mathbb Q}_{\alpha}= \{{\mathbb U}{_\alpha},\; {\mathbb
Y}_{\alpha}\}.$

Here, $\mathbb P_{\alpha}$= $\{\mathbb { U}{_\alpha},\; \mathbb X_{\alpha}\} $ is
defined by the axiumbilics  $\mathbb{U}_{\alpha}$ and the  field
of orthogonal tangent lines $\mathbb X_{\alpha}$, on  $\mathbb M\backslash  \mathbb{U}_{\alpha}$,
 on which the immersion is curved  along the large axis of the curvature ellipse.
The reason for the   name given to this object is that
for surfaces in $\mathbb R^3$,  $\mathbb P_{\alpha}$   reduces to the classical
{\it principal configuration} defined by the two principal curvature direction fields
$\{\mathbb X_{\alpha 1},\mathbb X_{\alpha 2}\}$,  \cite{gs1, gs2}.
Also, in  ${\mathbb Q}_{\alpha}= \{\mathbb { U}{_\alpha},\; \mathbb Y_{\alpha}\} $,
  $\mathbb Y_\alpha$ is the  field
of orthogonal tangent lines $\mathbb Y_\alpha$  on $\mathbb M\backslash \mathbb {U}_\alpha$,
on which the immersion is curved along the small axis of the curvature ellipse.
 For surfaces in $\mathbb R^3$    the curvature ellipse reduces to
  a segment and the crossing $\mathbb Y_\alpha$ splits into the two mean curvature line fields
    $\{\mathbb Y_{\alpha 1},\mathbb Y_{\alpha 2}\}$.
 In this case $\mathbb Q_{\alpha}$ reduces to the {\it mean configuration} defined by {\it umbilic points}
 and line fields along which the normal curvature is equal to the  Mean
 Curvature. That is the arithmetic mean of the principal
 curvatures.

 The global aspects of {\it arithmetic mean configurations} of surfaces immersed in $\mathbb R^3$
 have been studied  by Garcia and Sotomayor in \cite{m}.
  Examples of quadratic  ellipsoids such that all  arithmetic mean curvature lines are dense  were given.

Other mean curvature functions have been studied by Garcia and
Sotomayor in \cite{g, h, me}, unifying the arithmetic, geometric
and harmonic classical means of the principal curvatures.

The  global generic structure of the {\it axial principal} and {\it mean curvature}  lines,
 along which the second fundamental form points  in the direction of the large and the small
  axes of the Ellipse of
Curvature was developed by Garcia and Sotomayor, \cite{ax}.

\begin{figure}[htb]

 \begin{center}
  \includegraphics[angle=0, height=7cm, width=7cm]{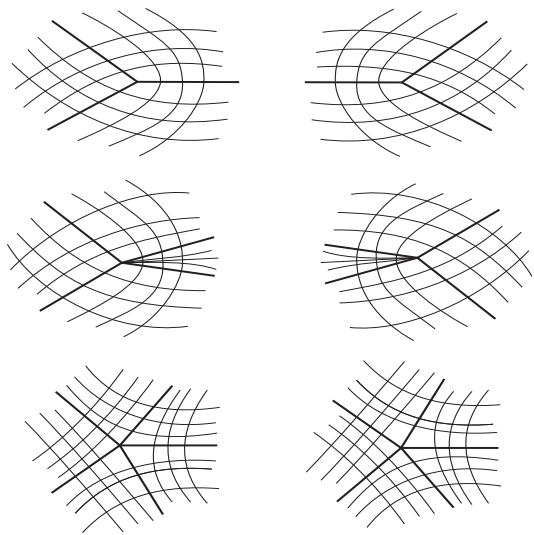}
  \caption{Axial configuration   near axiumbilic points.
  \label{fig:ax} }
 \end{center}
  \end{figure}

A   partial local attempt in this direction has been made by Gutierrez et
al.
in the
paper  \cite{ggtg}, where the structure around  the generic axiumbilic points
(for which the ellipse is a circle)  is established for surfaces  smoothly  immersed  in $\mathbb R^4$.
See Fig.  \ref{fig:ax} for an illustration of the three generic types: $E_3,\; E_4, \;E_5$,
established in \cite{ggtg} and also in \cite{ax}.

These points must regarded as the analogous to the Darbouxian umbilics:
{\rm D}$_1$,\; {\rm D}$_2$,\; {\rm D}$_3$, \cite{da} , \cite{gs1, gs2}. In both cases, the subindices
 refer to the number of {\it separatrices} approaching the singularity.

\newpage

\subsection{Differential equation for lines of axial curvature }\hskip 2cm
\medskip

Let $\alpha:\mathbb M^2\to \mathbb R^4$ be a $C^r,\;\; r\geq 4,$ immersion of an
 oriented smooth surface $\mathbb M$ into $\mathbb R^4$.
Let $N_1$ and $ N_2$ be a frame of    vector fields orthonormal to $\alpha$.
Assume that $(u,v)$ is a positive chart and that $\{\alpha_u, \alpha_v, N_1,N_2\}$
is  a positive frame.

In a chart $(u,v)$, the first fundamental form of $\alpha$ is
given by:

$I_\alpha= <D\alpha,D\alpha>= E du^2+2Fdudv+Gdv^2$, with

$E=<\alpha_u,\alpha_u>$, $F=<\alpha_u,\alpha_v>$, $G=<\alpha_v,\alpha_v>$

The second fundamental form is given by:
$$II_\alpha= II_{1,\alpha}N_1+ II_{2,\alpha}N_2$$ where,

 $II_{1,\alpha}= <N_1, D^2\alpha> = e_1 du^2 + 2f_1 dudv + g_1 dv^2$ and

 $II_{2,\alpha}= <N_2, D^2\alpha> = e_2 du^2 + 2f_2 dudv + g_2 dv^2$.

The normal curvature vector at a point $p$ in a tangent direction $v$ is given by:

 $ k_n= k_n(p,v) = II_\alpha(v,v)/I_\alpha(v,v).$

Denote   by $ T\mathbb M$ the tangent bundle of $\mathbb M$  and by $N\mathbb M$ the  normal
bundle of $\alpha$.
The image of the unitary circle of $T_p \mathbb M$ by
$k_n(p): T_p\mathbb M \to N_p\mathbb M$, being a quadratic,  map  is either
an ellipse,
a point or a segment.  In any case, to unify the notation, it will be
  refereed to as the {\it ellipse of curvature}
of $ \alpha$ and denoted by $\mathbb E_\alpha$.

The {\it mean curvature vector} ${\mathcal H}$ is defined by:
$${\mathcal H}=h_1 N_1+h_2 N_2=\frac{Eg_1+e_1G-2f_1F}{2(EG-F^2)}N_1+\frac{Eg_2+e_2G-2f_2F}{2(EG-F^2)}N_2.$$

Therefore,  the ellipse of curvature $\mathbb E_\alpha$ is given by the image of:

$$ k_n=(k_n-{\mathcal H})+ {\mathcal H}.$$

The tangent directions for which the  normal curvature  are  the axes, or vertices, of the ellipse of curvature
 $\mathbb E_\alpha$ are characterized  by the
following quartic form  given by the  Jacobian of the pair of forms
below, the first being quartic and the second quadratic:

$$Jac(||k_n-{\mathcal H}||^2, I_\alpha)=0.$$
\noindent where,
$$\aligned ||k_n-{\mathcal H}||^2 &=
 \big[\frac{e_1 du^2 + 2f_1 dudv + g_1 dv^2}{ E du^2+2Fdudv+Gdv^2}-\frac{(Eg_+e_1G-2f_1F)}{2(EG-F^2)}\big]^2\\
&+ \big[\frac{e_2 du^2 + 2f_2 dudv + g_2 dv^2}{ E du^2+2Fdudv+Gdv^2}-\frac{(Eg_2+e_2G-2f_2F)}{2(EG-F^2)}\big]^2. \endaligned $$

Expanding the equation above,  it follows that the
 differential equation  for the corresponding tangent directions, which defines the {\it   axial curvature lines},  is given by a quartic differential equation:

\begin{equation}\label{eq:dax}\aligned A &(u,v,du,dv)= [ a_0   G (EG-4F^2) + a_1 F(2F^2-EG) ] dv^4 \\
+&    [-8a_0EFG+ a_1E(4F^2-EG)]dv^3du \\
 +  & [  - 6a_0  GE^2+ 3a_1FE^2]dv^2 du^2
 +   a_1 E^3 dv du^3 + a_0 E^3 du^4  =0\endaligned\end{equation}

  \noindent where,
$$\aligned a_1 =& 4G(EG-4F^2)(e_1^2+e_2^2)
+32EFG(e_1f_1 +e_2f_2)  \\ +&  4E^3(g_1^2+g_2^2)
 -8E^2G(e_1g_1+e_2g_2)-   16E^2G(f_1^2+f_2^2)
\\
 a_0  = &4F(EG-2F^2)(e_1^2+e_2^2)  -4E(EG-4F^2)(e_1f_1+e_2f_2) \\
 -& 8E^2F(f_1^2+f_2^2)
 -4E^2F(e_1g_1+e_2g_2)+   4E^3(f_1g_1+f_2g_2).
\endaligned $$

\begin{remark} Suppose that  the surface $\mathbb M$ is contained into $\mathbb R^3$ with $e_2=f_2=g_2=0$. Then  the differential equation \eqref{eq:dax}
is the product of the differential equation of its {\it principal curvature lines}  and the differential equation of its {\it mean curvature  lines}, i.e.,
  the quartic differential equation \eqref{eq:dax} is  given by $$Jac(II_\alpha,I_\alpha).Jac(Jac(II_\alpha,I_\alpha),I_\alpha)=0. $$
\end{remark}

Let $\mathbb M^2$ be a compact, smooth and oriented  surface.
Call ${\mathcal M}^k$ the space of $C^k$ immersions of $\mathbb M^2$ into $\mathbb R^4$, endowed with the $C^k$ topology.

An immersion $\alpha\in {\mathcal M}^k$ is said to be {\it Principal Axial Stable} if it has a $C^k$, neighborhood ${\mathcal V}(\alpha)$, such that for any $\beta\in   \mathcal V(\alpha)$ there exist a homeomorphism $h:\mathbb M^2\to \mathbb M^2$ mapping $\mathbb U_\alpha$ onto  $\mathbb U_\beta$ and mapping the integral net of $\mathbb X_\alpha$ onto that of $\mathbb X_\beta$.
Analogous definition is given for {\it Mean Axial Stability}.

Sufficient conditions are provided to extend to the present setting  the Theorem on Structural Stability   for Principal Configurations due to  Gutierrez and Sotomayor \cite{gs1,   gsln, gs2}.
Consider the subsets $\mathcal P^k$ (resp. $\mathcal Q^k$) of immersions $\alpha$ defined by the following conditions:
\begin{itemize}
\item [a)] all axiumbilic points are of  types: $E_3$, $E_4$ or $E_5$;
\item[b)]  all principal (resp. mean) axial cycles  are hyperbolic;
\item[c)] the limit set of every axial line of curvature is contained in the set of axiumbilic points and principal (resp. mean) axial cycles of $\alpha$;
\item[d)] all axiumbilic separatrices are associated to a single axiumbilic point; this means that there are no connections or self connections of axiumbilic separatrices.
\end{itemize}

\begin{theorem} \label{th:ax}\cite{ax} Let $k\geq 5$.  The following holds:
\begin{itemize}
\item[i)] The subsets ${\mathcal P}^k$ and $ {\mathcal Q}^k$ are  open in  ${\mathcal M}^k$;
\item[ii)] Every $\alpha\in {\mathcal P}^k$ is  Principal Axial Stable;
\item[iii)] Every $\alpha\in{ \mathcal Q}^k$ is  Mean Axial Stable.
\end{itemize}
\end{theorem}

\begin{remark}
Several other extensions of principal lines of surfaces immersed in $\mathbb R^4$ have been considered.
To have an idea of these developments the reader is addressed to \cite{gms}, \cite{sanchez-garcia}, \cite{gugui}, \cite{mello}, \cite{sanchez}.
\end{remark}

 \section{ Principal Configurations on Immersed Hypersurfaces in $\mathbb R^4$}

Let $ \mathbb M^m$  be a $C^k, k\ge 4$, compact and oriented,
$m-$~dimensional manifold.
An immersion $\alpha$ of $ \mathbb M^m$ into $\re^{m+1}$ is a map such that
$D{\alpha}_p : T\mathbb M^m_p  \to \re^{m+1} $ is one to one, for every $p\in
\mathbb M^m$.
 Denote by ${\mathcal J}^{k}  = {\mathcal J}^{k} (\mathbb M^m, \re^{m+1}) $ the
set
of $C^{k}$ - immersions of $ \mathbb M^m $ into $\re^{m+1}$. When endowed with
the $C^s-$ topology of Whitney, $s\le k, $ this set is denoted by $
{\mathcal J} ^{k,s}$. Associated to every $\alpha \in {\mathcal J}^k$ is defined the
normal map $N_{\alpha}: \mathbb M^m \to S^m$ :
$$ N_{\alpha} =  ({{\alpha}_1\we  \ldots \we{\alpha}_m}) /{\mid
{\alpha}_1\we  \ldots \we{\alpha}_m\mid },$$
\noindent  where $({u_1},\ldots,{u_m}): (M , p) \to (\mathbb R^m , 0) $
 is a positive chart of $ \mathbb M^m$  around p, $\wedge$ denotes the exterior
product of vectors
in $\mathbb R^{m+1}$  determined by a once for all fixed orientation of
$\mathbb R^{m+1}$,
 ${\alpha}_1 = \prf{\alpha}{u_1}$, $\ldots$, $ {\alpha}_m = \prf{\alpha}{u_m} $ and
$\mid$
\hskip .3cm $\mid = {< , >}^{\frac 12}$ is the Euclidean norm in $\mathbb R^{m+1}$.
Clearly, $ N_{\alpha}$ is well defined  and of  class $C^{k-1}$ in
$\mathbb M^m$.

Since $D N_{\alpha}(p)$ has its image contained in that of $D\alpha
(p),$ the endomorphism
${\omega}_{\alpha} : T\mathbb M^m \to T\mathbb M^m $ is well
defined by
$$ D\alpha. {\omega}_{\alpha} = D N_{\alpha}. $$
It is well know that ${\omega}_{\alpha}$ is  a self  adjoint endomorphism,
when $ T\mathbb M^m $ is endowed with the metric ${< , >}_{\alpha} $ induced by
$\alpha$ from the metric in $\re^{m+1}$.

The opposite values of the eigenvalues of ${\omega}_{\alpha}$ are
called {\it principal curvatures} of  $\alpha$ and will be denoted by $k_1
\le \ldots \le k_m$.
 The eigenspaces associated to the principal
curvatures define m $\;\;C^{k-2}$ line fields \lp i, $(i=1,\ldots,m)$
 mutually orthogonal in $T\mathbb M^m$ (with the metric $ { < ,
>}_{\alpha}$), called {\it principal line fields} of $\alpha$. They are
characterized by Rodrigues' equations \cite{Sp}, \cite{St}.
{\small $$ {{\mathcal L}}_1(\alpha)  = \{ v\in T\mathbb M^m : {\omega}_{\alpha}v +
k_1 v = 0\},  \ldots,
{{\mathcal L}}_m(\alpha) = \{ v\in TM : {\omega}_{\alpha}v + k_m v =
0\}.      $$}
The integral curves of \lp i, $(i=1, \ldots, m) $  outside their
singular set,  are
called {\it lines of principal curvature}. The family of such
curves i.e. the integral foliation of \lp i
will be denoted by \fp i and are called the
{\it  principal foliations} of $\alpha$.

\subsection{Curvature lines near Darbouxian partially umbilic curves}\hskip 2cm
\medskip

In the three dimensional case, there are three principal foliations \fp
i which are mutually orthogonal. Here two kind of singularities of the principal
line fields \lp i  $(i=1, 2, 3)$ can appear. Define the sets,
${\mathcal U}(\alpha)=\{ p\in \mathbb M^3: k_1(p) = k_2(p) = k_3(p)\}$,
${\mathcal P}_{12}(\alpha)=\{ p\in \mathbb M^3: k_1(p) = k_2(p) \ne k_3(p)\}$,

${\mathcal P}_{23}(\alpha)=\{ p\in \mathbb M^3: k_1(p) \ne k_2(p) = k_3(p)\}$ and
${\mathcal P}(\alpha)={\mathcal P}_{12}(\alpha)\cup {\mathcal P}_{23}(\alpha)$.

The sets ${\mathcal U}(\alpha)$,  ${\mathcal P}(\alpha)$ are called,  respectively,{\it
umbilic set } and  {\it partially umbilic set } of the immersion $\alpha$.

Generically, for an open and dense set of immersions in the space
${\mathcal J}^{k,s}$,
${\mathcal U}(\alpha)= \emptyset$ and  ${\mathcal P}(\alpha)$ is either, a submanifold of
codimension two or the empty set.

A connected component of ${\mathcal  S}(\alpha)$ is called a {\it partially umbilic curve}.

The study of the principal foliations near $\mathcal S (\alpha)$ were
carried
out in  \cite{ga-tese, ga1, ga2}, where the local model of the  asymptotic behavior of
 lines of principal curvature was   analyzed in the generic
case.

In order to state the results  the following definition  will be introduced.

\begin{definition} \label{def:1r4} Let $ p\in {\mathbb M}^3$ be  a partially
umbilic point such that $ k_1(p) \ne k_2(p) =k_3(p) =k(p)$.

\noindent Let $\uu: {\mathbb M}^3 \rightarrow \mathbb R^3 $ be a local chart
and  $ R: \mathbb R^4 \to \mathbb R^4 $ be an  isometry such that:

 $ (R\circ \alpha)\uu = (u_1, u_2, u_3, h\uu)=u_1 e_1+u_2 e_2+u_3 e_3 + h\uu e_4$,
 where $\{e_1,e_2,e_3,e_4\}$ is the canonical basis of $\mathbb R^4$ and
{\small
$$\aligned  h\uu &= \dfrac 12k_1\ue 21 + \dfrac 12k (\ue 22 +
\ue 23) + \dfrac  16 a_1\ue 31 + \dfrac 12 a_2 \ue 21 \ui 2 + \dfrac
12 a_3 \ue 21\ui 3 \\&+ \dfrac 12 a_4 \ui 1\ue 22 +  \dfrac
12 a_5 \ui 1\ue
23 + a_6 \ui 1 \ui 2 \ui3 \\
 &+  \dfrac 16 a\ue 32 + \dfrac 12 b\ui 2\ue 23
+ \dfrac 12 c\ue 33   + h.o.t.
\endaligned $$
}

The point $p$  is called a {\it Darbouxian  partially umbilic point} of  type {\rm D}$_i$ if the $p$ verifies the transversality condition
$T)\quad   b (b - a) \ne 0 $   and the condition  {\rm D}$_i$ below holds
 {\small
$$ \aligned \text{D}_1)\qquad \frac ab &> (\frac{c}{2b})^2 + 2 \\
\text{D}_2)\qquad 1 &< \frac ab  < (\frac{c}{2b})^2 + 2, \qquad a \ne 2b \\
\text{D}_3)\qquad \frac ab  &< 1 \endaligned$$ } \end{definition}

 We observe that the term $d u_2^2 u_3$ was eliminated by an appropriated rotation.

A {\em  partially umbilic point} which satisfy the  condition $T$   belongs to a regular curve formed by  partially umbilic points and that this curve is transversal to the geodesic surface tangent to the {\it umbilic plane}, which is the plane spanned by the eigenvectors corresponding to the multiple eigenvalues.

A  regular arc of Darbouxian partially umbilic points  {\rm D}$_i$  will be  called {\it Darbouxian partially umbilic curve } {\rm D}$_i$ .

\begin{remark}
\begin{itemize}

\item[i)] Along a regular  connected component of  ${\mathcal P}_\alpha$ it is expected
that
can occur transitions of the types {\rm D}$_i$, $i = 1, 2, 3.$
 The simpler
transitions are those of  types {\rm D}$_{1}$-{\rm D}$_{2}$ and {\rm D}$_{2}$-{\rm D}$_{3}$. See \cite{ga-tese} and \cite{ga1}.

\item[ii)]  The conditions that define the types {\rm D}$_i$   are
independent of   coordinates. Also, these conditions are
closely related to those that define the {\it Darbouxian   umbilic points} in the two dimensional case, see \cite{bf}, \cite{da} and \cite{sell, gs1, gs-5, gs2}, as well as subsection \ref{ssec:cod1}.

\end{itemize}
 \end{remark}

The main result  announced in \cite{ga1} and proved in  \cite{ga-tese,ga2} is  the following theorem.

\begin{figure}[htbp]
\begin{center}
\includegraphics[angle=0, height= 3.0cm, width=10.0cm ]{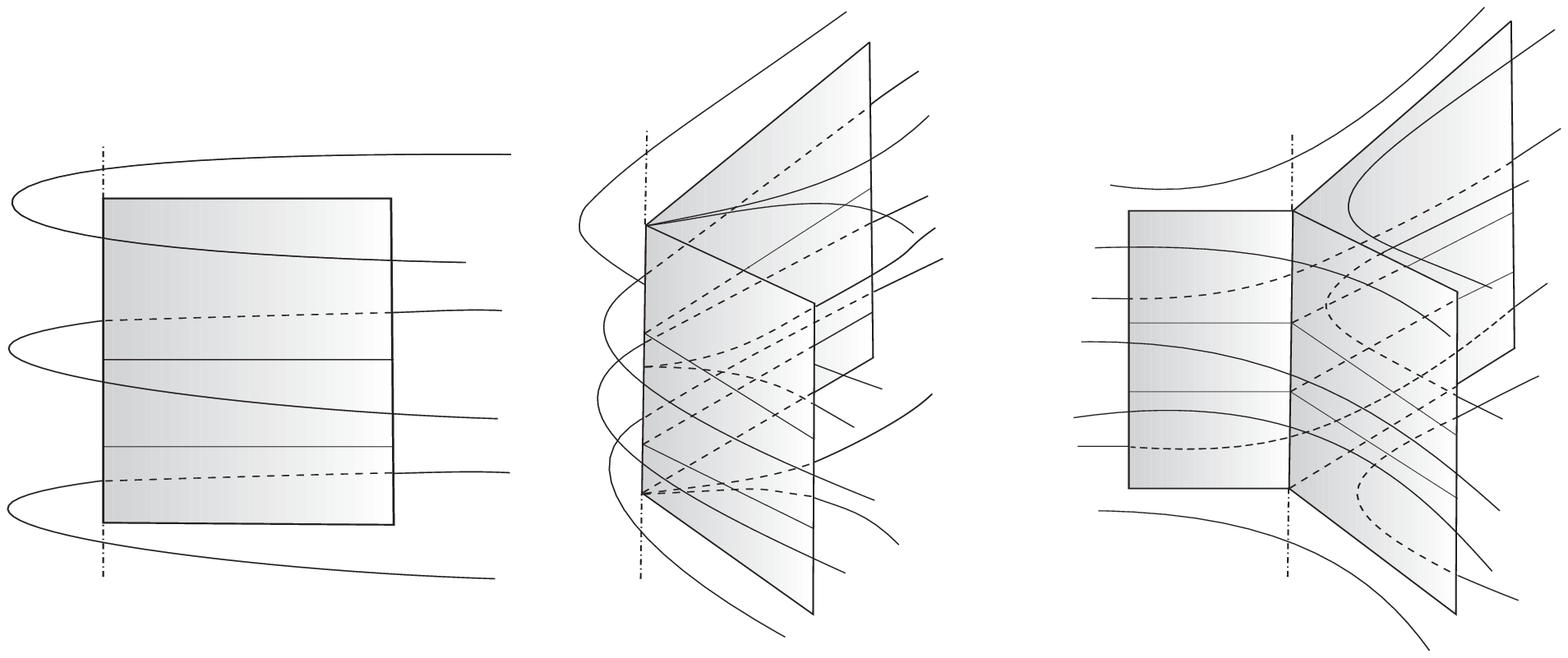}
\end{center}
\caption{ Behavior of  \fp 2  in
the neighborhood of a \lise\; {\rm D}$_1, $  {\rm D}$_2$ and {\rm D}$_3.$
\label{fig:t-d123}
  }
\end{figure}

\begin{figure}[htbp]
\begin{center}
\includegraphics[angle=0, height= 3.5cm, width=10.0cm ]{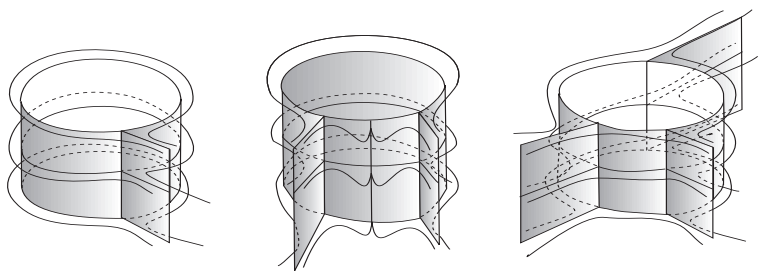}
\end{center}
\caption{ Behavior of \fp 2 near the projective line (blowing-up).
\label{fig:d123p}}
\end{figure}

 \begin{theorem}  \label{th:pul} \cite{ga-tese, ga1,ga2} Let $\alpha \in {\mathcal J}^k({\mathbb M}^3, \mathbb R^4), \;\; k\geq 4,$
and $p$ be a Darbouxian partially umbilic point.
Let $p \in c$, where $ c$ is a  Darbouxian partially
umbilic curve and $V_c$ a tubular neighborhood of
 $c$. Then it follows that:
 \begin{itemize}
 \item [i)] If $c$ is a
Darbouxian  partially umbilic curve {\rm D}$_1$, then there exists  an unique
invariant surface $W_c$ (umbilic separatrix)
of class  $C^{k-2},\;$ fibred over $c$ and whose fibers are leaves
of \fp 2 and the boundary of $W_c$ is $c$.
The set $V_c\backslash W_c$ is a  hyperbolic region  of  \fp 2.

 \item [ii)]  If $c$ is a
Darbouxian  partially umbilic curve  {\rm D}$_2$, then
there exist  two surfaces as above and exactly one parabolic region and
one hyperbolic region of \fp 2.
 \item [iii)]  If $c$ is a
Darbouxian  partially umbilic curve {\rm D}$_3$, then there exist
three  surfaces  as above and exactly three  hyperbolic regions of \fp2.

 \item [iv)]   The same happens for the foliation  \fp 3 which is
orthogonal to \fp 2 and singular in the curve $c$. Moreover,  the invariant  surfaces of  \fp 3, in each case {\rm D}$_i$, are tangent to the invariant surfaces of \fp 2 along $c$.

 \item [v)]  Fig.  \ref{fig:t-d123}   shows the behavior of \fp 2 in the neighborhood of a Darbouxian partially umbilic curve. The foliation \fp 3 has a similar behavior.  Fig.  \ref{fig:d123p} shows the behavior of \fp 2 near the projective line (blowing-up).
     \end{itemize}
\end{theorem}

\vspace{-.3cm}

    \begin{remark} \begin{itemize}
\item[i)] The union of the invariant surfaces of \fp 2 and \fp 3 which have the same tangent plane at the curve $c$ is of class $C^{k-2}$.

\item[ii)] The distribution of planes defined by
\lp i, i.e. the distribution that  has \lp i as a normal vector, is
not
integrable in general, and so, the situation is strictly tridimensional.
\end{itemize}
\end{remark}

\begin{figure}[hb]
\begin{center}
\includegraphics[angle=0, height= 3.2cm  ]{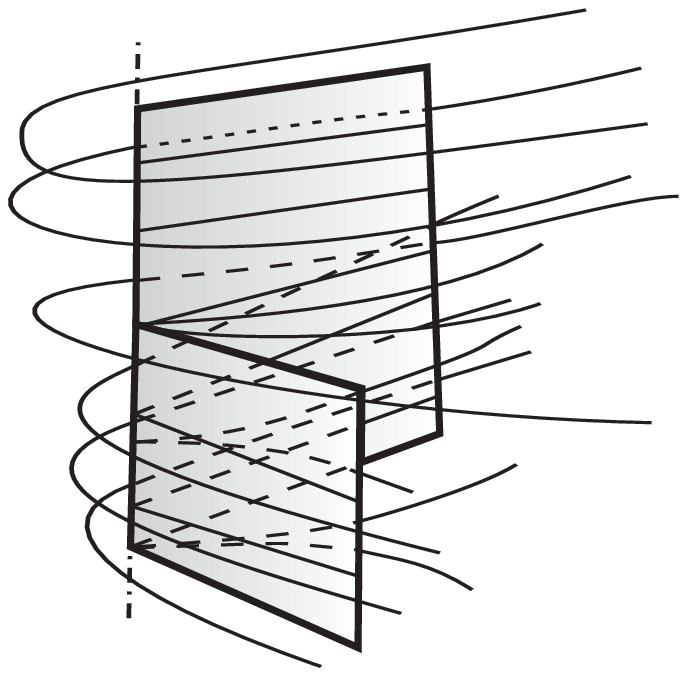}
\includegraphics[angle=0, height= 3.8cm, width=10cm ]{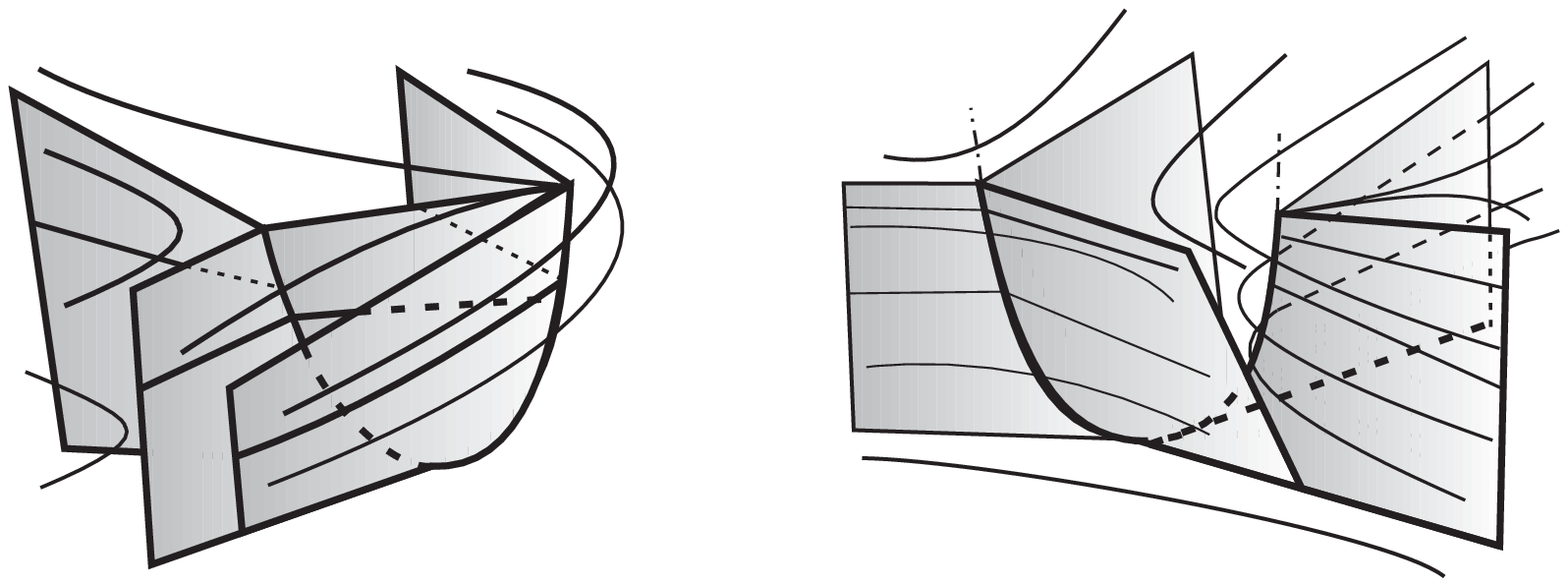}
\end{center}
\caption{ Behavior of \fp 2 near a transition {\rm D}$_1$-{\rm D}$_2$ and  of \fp 2 and \fp 3 near a transition {\rm D}$_2$-{\rm D}$_3.$
\label{fig:d12d23}}
\end{figure}

In  Fig.    \ref{fig:d12d23}   it is shown the generic behavior of principal foliation \fp 2 and \fp 3 near the  transitions of arcs {\rm D}$_1$-{\rm D}$_2$ and {\rm D}$_2$-{\rm D}$_3$. See \cite{ga-tese, ga1}.

\subsection{Curvature lines near hyperbolic principal cycles}\hskip 2cm
\medskip

Next it will be described the behavior of principal foliations near principal cycles, i.e., closed principal lines of \fp i.

Let  $c$ be principal cycle of the principal foliation \fp 1.

With respect to the positive orthonormal frame $\{ e_1, e_2 , e_3 ,N \},$   Darboux's equations for the curve
$\alpha\circ c$ are given by the following system, see \cite{Sp}.
$$\aligned \left(\begin{matrix} e_1^\prime\\ e_2^\prime \\ e_3^\prime \\
N^\prime\end{matrix}\right)
&=\left(\begin{matrix} 0 & \oa 12 & \oa 13 & k_1 \\ -\oa 12 & 0 & \oa 23 & 0 \\
-\oa 13 & -\oa 23 & 0 & 0 \\ - k_1 &  0 & 0 & 0 \end{matrix}\right) \left(\begin{matrix}  e_1
\\ e_2 \\ e_3 \\ N \end{matrix}\right)\endaligned  $$
where $\oa ij = < \nabla_{e_1}e_i , e_j >. $

In order to study the behavior of  \fp 1 in the neighborhood
of $c$ we will study the    Poincar\'e map
$\Pi$ associated to \fp 1 whose definition is reviewed as follows.

Let  $\uu$  be the system of  coordinates given by lemma 1.
Consider  in these  coordinates the transversal sections
$\quad\{ \ui 1 = 0 \}$ and
$\{ \ui 1 = L \}\quad$ and define the map  $ \Pi :\{ \ui 1
= 0 \} \rightarrow \{ \ui 1 = L \}$ in the following way.

\noindent Suppose $ c $ be oriented by the parametrization $\ui
1$ and let $\uu $ the solution of  the differential equation that
defines the  principal line field \lp 1,
with initial condition
$\uu  (0 ,\ue 02 ,\ue 03) =  (0 ,\ue 02 ,\ue 03) $.
So the Poincar\'e map is defined by

  $$\Pi (\ue 02 ,\ue 03) = (\ui 2\scriptstyle (L  ,\ue 02 ,\ue 03)\displaystyle , \ui 3 \scriptstyle
 (L ,\ue 02 ,\ue 03)\displaystyle).$$

The principal cycle $c$ is called hyperbolic if the eigenvalues of $\Pi^\prime (0)$ are disjoint from the unitary circle, \cite{soto}.

\begin{proposition} In the conditions above
 we have that the derivative of the Poincar\'e map $\Pi$
is given by $\Pi^\prime (0) = U (L)$, where  $U $ is the solution
of the following differential equation:

$$\aligned    \left\{ \begin{matrix}    U^\prime &= A U \qquad \qquad  \qquad    \qquad\quad \\
U (0)&= I,\,   A (\ui 1) = A (\ui 1 + L), \end{matrix} \right.
  \;  A(\ui 1) &=
\left(\begin{matrix}
\dfrac {\scriptstyle -k_2^\prime\displaystyle}{\scriptstyle (k_2 - k_1)
\displaystyle}
& \dfrac{\scriptstyle \oa 23(k_3 - k_1)\displaystyle}{\scriptstyle (k_2
- k_1)\displaystyle}\cr
\dfrac{\scriptstyle-\oa 23 (k_2 - k_1)\displaystyle }{\scriptstyle (k_3
- k_1)\displaystyle} & \dfrac{\scriptstyle  -
k_3^\prime\displaystyle}{\scriptstyle (k_3
- k_1)\displaystyle }\end{matrix}\right). \endaligned
  $$

\end{proposition}

 \begin{theorem}\cite{gagag} Let $c$ be a  principal cycle of  ${\mathcal F}_1 (\alpha).$
Suppose that  the principal curvature $k_1$ is not constant along  $c$.
Then, given $\epsilon > 0 $, there exists an immersion $\tilde\alpha \in
{\mathcal J}^{\infty,r} (\mathbb \mathbb M^3 ,\mathbb R^4
) $ such that $||\alpha-\tilde\alpha||<\epsilon$  and $c$ is a
hyperbolic  principal cycle of  ${\mathcal F}_1 (\tilde\alpha).$
\end{theorem}

\begin{remark} Global aspects of principal lines of immersed hypersurfaces in $\mathbb R^4$ was studied in \cite{ga-tese}.\end{remark}

\section{Concluding Comments
}\label{sec:cc}

In this work an effort has been made to present  most of  the developments  addressed
to improve the  local and global understanding of
 the structure of principal curvature lines
on surfaces and related geometric objects. The emphasis has been put on those developments    derived from the assimilation   of ideas coming from the {\it QTDE}
and Dynamical Systems
into the classical knowledge on  the subject, as  presented in prestigious treatises such as Darboux \cite{da},  Eisenhart \cite{ei},  Struik \cite{St},  Hopf \cite{hopf},
Spivak \cite{Sp}.

The starting point for the results surveyed here can be found in the papers of Gutierrez
and Sotomayor \cite{gs1, gsln, gs2}, as suggested in the historical essay contained
in sections \ref{sec:2} to \ref{sec:gs}.

The authors acknowledge the influence  they received from the well established
 theories of Structural Stability
and Bifurcations which developed from  the inspiring classical works of
 Andronov,  Pontrjagin, Leontovich \cite{ap} and Peixoto \cite{mp, ppr}.
  Also the results on bifurcations of principal configurations  outlined in \cite{gs-5} and further elaborated along this work are motivated in Sotomayor \cite {soihes}.

The vitality of the  {\it QTDE} and Dynamical Systems, with  their remarkable present day achievements, may lead to the belief  that the possibilities for di\-rections of future research
in connection with the differential equations of lines of curvature and other equations of
Classical Geometry   are too wide and undefined and that the source of problems in the subject
  consists mainly in establishing an analogy with  one in the above mentioned fields.

While  this may partially true in the present work, History shows us that the consideration
 of problems derived from  purely geometrical  sources  and from other fields
 such as Control Theory, Elasticity, Image Recognition and  Geometric Optics,
  have also a crucial role to play in determining the directions for relevant research in our subject.
In fact, at the very beginning, the works of Monge and
Dupin and, in relatively recent times, also   the famous Carath\'eodory Conjecture
\cite{gagu}, \cite{gmb}, \cite{r}, \cite{ha1,ha2,ha3}, \cite{iva},  \cite{meso},
\cite{taba, taba2}, \cite{sx3}, \cite{sx1,sx2},  \cite{xa},    represent geometric sources
of research
directions leading to clarify the structure of lines of curvature and their umbilic singularities.

 \vskip 0.2cm
  \noindent{\bf Acknowledgement.\,}  The authors are grateful
 to   L. F. Mello and C. Gutierrez   for   helpful
 remarks.  Sotomayor thanks the hospitality of the Mathe\-ma\-tics  Department at Brown University,
 where part of this work was done.

\vskip .5cm
\author{\noindent Jorge Sotomayor\\Instituto de Matem\'{a}tica e Estat\'{\i}stica,\\Universidade de S\~{a}o Paulo,
\\Rua do Mat\~{a}o 1010, Cidade Universit\'{a}ria, \\CEP 05508-090, S\~{a}o Paulo, S.P., Brazil \\
\\ Ronaldo Garcia\\Instituto de Matem\'{a}tica e Estat\'{\i}stica,\\
Universidade Federal de Goi\'as,\\CEP 74001-970,
Caixa Postal 131,\\Goi\^ania, GO, Brazil}

\end{document}